\documentclass[11pt]{article}
\usepackage{amsmath}
\usepackage{amssymb}
\usepackage{amscd}
\usepackage{xspace}
\usepackage{verbatim}
\newcommand{\mysection}[1]{
\section{#1}\setcounter{equation}{0}}
\date{}
\begin{document}
\begin{center}{\bf \Large BOUNDARY VALUE PROBLEMS WITH MEASURES FOR \\[3mm]

\noindent ELLIPTIC EQUATIONS WITH SINGULAR POTENTIALS}
\footnote{Both authors are sponsored by the ECOS-Sud program C08E04. The second author is partially supported by Fondecyt 1070125}
 \end{center}
 %% Title to be inserted !!!!!!!!!!
\begin{center}{\bf  Laurent V\'eron}\\
{\small Laboratoire de Math\'ematiques et Physique Th\'eorique}\\
 {\small  Universit\'e Fran\c cois-Rabelais, Tours,  FRANCE} \\[2mm]

 {\bf  Cecilia Yarur}\\
 {\small Departamento de Matem\'aticas y Ciencia de la Computaci\'on} \\
 {\small  Universidad de Santiago de Chile, Santiago, CHILE }\\[2mm]
with an appendix by {\bf  Alano Ancona}\\
 \end{center}

%\maketitle
%%

%%
%\maketitle

%% FONT commands
\newcommand{\txt}[1]{\;\text{ #1 }\;}%% Used in math only
\newcommand{\tbf}{\textbf}%% Bold face. Usage: \tbf{...}
\newcommand{\tit}{\textit}%% Italic
\newcommand{\tsc}{\textsc}%% Small caps
\newcommand{\trm}{\textrm}
\newcommand{\mbf}{\mathbf}%% Math bold
\newcommand{\mrm}{\mathrm}%% Math Roman
\newcommand{\bsym}{\boldsymbol}%% Bold math symbol
%%Macros for changing font size in math.
\newcommand{\scs}{\scriptstyle}%% as in subscript
\newcommand{\sss}{\scriptscriptstyle}%% as in sub-subscript
\newcommand{\txts}{\textstyle}
\newcommand{\dsps}{\displaystyle}
%%Macros for changing font size in text.
\newcommand{\fnz}{\footnotesize}
\newcommand{\scz}{\scriptsize}
%%\tiny<\scz<\fsz<\small<\large<\Large<\huge<\Huge
%%%%%%%%%%%%
%%%%%%%%%%%%
%% EQUATION commands
\newcommand{\be}{
\begin{equation}
}
\newcommand{\bel}[1]{
\begin{equation}
\label{#1}}
\newcommand{\ee}{
\end{equation}
}%% This macro does not work with amstex.
\newcommand{\eqnl}[2]{
\begin{equation}
\label{#1}{#2}
\end{equation}
}%%use not advisable; confusing
%%%%%%%%%%%%%%%
%% Unnumbered THEOREM env.
%% New env. to be used for unnumbered theorem, lemma etc. (but with specified name)
\newtheorem{subn}{\name}
\newcommand{\bsn}[1]{\def\name{#1}
\begin{subn}}
\newcommand{\esn}{
\end{subn}}
%%%%%%%%%%%%%%
%% NUMBERED THEOREM env.
%% Environments: theorem, lemma, corollary defintion and related commands,
%% designed to provide consecutive numbering of these forms.
\newtheorem{sub}{\name}[section]
\newcommand{\dn}[1]{\def\name{#1}}   %used in conjuction with sub or subn.
\newcommand{\bs}{
\begin{sub}}
\newcommand{\es}{
\end{sub}}
\newcommand{\bsl}[1]{
\begin{sub}\label{#1}}
%% the above must be preceeded by \dn (name definition),
%% however this is superceded by the list of commands bth etc.  below.
%%%%%%%%%%%%
%% NUMBERED THEOREM env. (cont.)
%% List of commands derived from 'sub' env. for theorem, lemma etc.
%% designed to provide consecutive numbering of these forms.
\newcommand{\bth}[1]{\def\name{Theorem}
\begin{sub}\label{t:#1}}
\newcommand{\blemma}[1]{\def\name{Lemma}
\begin{sub}\label{l:#1}}
\newcommand{\bcor}[1]{\def\name{Corollary}
\begin{sub}\label{c:#1}}
\newcommand{\bdef}[1]{\def\name{Definition}
\begin{sub}\label{d:#1}}
\newcommand{\bprop}[1]{\def\name{Proposition}
\begin{sub}\label{p:#1}}
%%%%%%%%%%%%%%%%%%%%%%%%%%%%%%%%%%
%% RERERENCE commands.
%% \newcommand{\R}[1]{(\ref{#1})}
\newcommand{\R}{\eqref}
\newcommand{\rth}[1]{Theorem~\ref{t:#1}}
\newcommand{\rlemma}[1]{Lemma~\ref{l:#1}}
\newcommand{\rcor}[1]{Corollary~\ref{c:#1}}
\newcommand{\rdef}[1]{Definition~\ref{d:#1}}
\newcommand{\rprop}[1]{Proposition~\ref{p:#1}} 
%%%%%%%%%%%
%% ARRAY commands.
\newcommand{\BA}{
\begin{array}}
\newcommand{\EA}{
\end{array}}
\newcommand{\BAN}{\renewcommand{\arraystretch}{1.2}
\setlength{\arraycolsep}{2pt}
\begin{array}}
\newcommand{\BAV}[2]{\renewcommand{\arraystretch}{#1}
\setlength{\arraycolsep}{#2}
\begin{array}}
%Note: The first variable gives the amount of stretching: (#1) x default.
%For instance #1=1.2 means a 20% stretching. The second variable should be
%written for instance in the form  4pt ; here the default is 5pt
%\newcommand{\EAN}{\end{array}\setlength{\arraycolsep}{5pt}}
\newcommand{\BSA}{
\begin{subarray}}
\newcommand{\ESA}{\end{subarray}}
%Note: These are used in subscripts as well as superscripts. They work essentially
%% like 'array'.
\newcommand{\BAL}{\begin{aligned}}
\newcommand{\EAL}{\end{aligned}}
\newcommand{\BALG}{\begin{alignat}}
\newcommand{\EALG}{\end{alignat}}%% the abbrev. does not work with latex2e
\newcommand{\BALGN}{\begin{alignat*}}
\newcommand{\EALGN}{\end{alignat*}}%% the abbrev. does not work with latex2e
%% The 'aligned' environment must be placed inside an 'equation' env.
%% in the same way as the array.
%% One could use also the 'align' env. or the 'alignat' env.
%% However in this case each line is numbered, unless '\notag' is used.
%% The 'alignat'
%% has a slightly different format (the number of columns must be specified in advance)
%% but it has the advantage that the distance between columns is at our disposition.
%% (The default would be zero distance.) Using 'alignat*' we can have the advantages
%% of alignat plus the situation where separate lines are not numbered.
%% However in this case there is no numbering at all (unless we provide a tag).
%%%%%%%%%%
%% PROOF, REMARK etc.
\newcommand{\note}[1]{\textit{#1.}\hspace{2mm}}
\newcommand{\Proof}{\note{Proof}}
\newcommand{\qeda}{\hspace{10mm}\hfill $\square$}
\newcommand{\qed}{\\
${}$ \hfill $\square$}
\newcommand{\Remark}{\note{Remark}}
%%%%%%%% Style command.
\newcommand{\modin}{$\,$\\
[-4mm] \indent}
%% To be used after \mysection in order to start new line with \indent.
%%%%%%%%%%%%
%% MATHEMATICAL symbols
\newcommand{\forevery}{\quad \forall}
\newcommand{\set}[1]{\{#1\}}
\newcommand{\setdef}[2]{\{\,#1:\,#2\,\}}
\newcommand{\setm}[2]{\{\,#1\mid #2\,\}}
%% Arrows
\newcommand{\lra}{\longrightarrow}
\newcommand{\sgn}{\rm{sgn}}
\newcommand{\lla}{\longleftarrow}
\newcommand{\llra}{\longleftrightarrow}
\newcommand{\Lra}{\Longrightarrow}
\newcommand{\Lla}{\Longleftarrow}
\newcommand{\Llra}{\Longleftrightarrow}
\newcommand{\warrow}{\rightharpoonup}
%% Brackets, delimiters
\newcommand{
\paran}[1]{\left (#1 \right )}%% adjustable parantheses
\newcommand{\sqbr}[1]{\left [#1 \right ]}%% adjustable square brackets
\newcommand{\curlybr}[1]{\left \{#1 \right \}}%% adjustable curly brackets
\newcommand{\abs}[1]{\left |#1\right |}%% adjustable vertical delimiters
\newcommand{\norm}[1]{\left \|#1\right \|}%% adjustable norm
\newcommand{
\paranb}[1]{\big (#1 \big )}%% non-adjustable parantheses (big)
\newcommand{\lsqbrb}[1]{\big [#1 \big ]}%% non-adjustable square brackets (big)
\newcommand{\lcurlybrb}[1]{\big \{#1 \big \}}%% non-adjustable curly brackets (big)
\newcommand{\absb}[1]{\big |#1\big |}%% non-adjustable vertical delimiters (big)
\newcommand{\normb}[1]{\big \|#1\big \|}%% non-adjustable norm (big)
\newcommand{
\paranB}[1]{\Big (#1 \Big )}%% non-adjustable parantheses (Big)
\newcommand{\absB}[1]{\Big |#1\Big |}%% non-adjustable vertical delimiters (Big)
\newcommand{\normB}[1]{\Big \|#1\Big \|}%% non-adjustable norm (Big)

%%%%%%%%%%%%%%%%%
%% Adjustable parantheses etc. in a different DEFINITION format.
%\def\adp(#1){\left (#1 \right )}%% adjustable parantheses
%\def\adsb(#1){\left [#1\right ]}%% adjustable square brackets
%\def\adcb(#1){\left \{#1\right \}}%% adjustable curly brackets
%\def\abs|#1|{\left |#1\right |}%% adjustable vertical delimiters
%%%%%%%%%%%%%%%%
%% More mathematical symbols
\newcommand{\thkl}{\rule[-.5mm]{.3mm}{3mm}}
\newcommand{\thknorm}[1]{\thkl #1 \thkl\,}
\newcommand{\trinorm}[1]{|\!|\!| #1 |\!|\!|\,}
\newcommand{\bang}[1]{\langle #1 \rangle}%% angle bracket
\def\angb<#1>{\langle #1 \rangle}%% angle bracket
%% The two last lines yield the same result.
%% The second is used as follows: \angb<a,b>
\newcommand{\vstrut}[1]{\rule{0mm}{#1}}
\newcommand{\rec}[1]{\frac{1}{#1}}
%% OPERATOR names.
%% OPERATOR names.
\newcommand{\opname}[1]{\mbox{\rm #1}\,}
\newcommand{\supp}{\opname{supp}}
\newcommand{\dist}{\opname{dist}}
\newcommand{\myfrac}[2]{{\displaystyle \frac{#1}{#2} }}
\newcommand{\myint}[2]{{\displaystyle \int_{#1}^{#2}}}
\newcommand{\mysum}[2]{{\displaystyle \sum_{#1}^{#2}}}
\newcommand {\dint}{{\displaystyle \int\!\!\int}}%%%%%%%%%%
%%%%%%% SPACE commands
\newcommand{\q}{\quad}
\newcommand{\qq}{\qquad}
\newcommand{\hsp}[1]{\hspace{#1mm}}
\newcommand{\vsp}[1]{\vspace{#1mm}}
%%%%%%%%%%%
%% ABREVIATIONS
\newcommand{\ity}{\infty}
\newcommand{\prt}{
\partial}
\newcommand{\sms}{\setminus}
\newcommand{\ems}{\emptyset}
\newcommand{\ti}{\times}
\newcommand{\pr}{^\prime}
\newcommand{\ppr}{^{\prime\prime}}
\newcommand{\tl}{\tilde}
\newcommand{\sbs}{\subset}
\newcommand{\sbeq}{\subseteq}
\newcommand{\nind}{\noindent}
\newcommand{\ind}{\indent}
\newcommand{\ovl}{\overline}
\newcommand{\unl}{\underline}
\newcommand{\nin}{\not\in}
\newcommand{\pfrac}[2]{\genfrac{(}{)}{}{}{#1}{#2}}% frac with parantheses.
%%%%%%%%%%%
%%%%%%%%%%%%%

%%Macros for Greek letters.
\def\ga{\alpha}     \def\gb{\beta}       \def\gg{\gamma}
\def\gc{\chi}       \def\gd{\delta}      \def\ge{\epsilon}
\def\gth{\theta}                         \def\vge{\varepsilon}
\def\gf{\phi}       \def\vgf{\varphi}    \def\gh{\eta}
\def\gi{\iota}      \def\gk{\kappa}      \def\gl{\lambda}
\def\gm{\mu}        \def\gn{\nu}         \def\gp{\pi}
\def\vgp{\varpi}    \def\gr{\rho}        \def\vgr{\varrho}
\def\gs{\sigma}     \def\vgs{\varsigma}  \def\gt{\tau}
\def\gu{\upsilon}   \def\gv{\vartheta}   \def\gw{\omega}
\def\gx{\xi}        \def\gy{\psi}        \def\gz{\zeta}
\def\Gg{\Gamma}     \def\Gd{\Delta}      \def\Gf{\Phi}
\def\Gth{\Theta}
\def\Gl{\Lambda}    \def\Gs{\Sigma}      \def\Gp{\Pi}
\def\Gw{\Omega}     \def\Gx{\Xi}         \def\Gy{\Psi}

%%Macros for calligraphic letters.
\def\CS{{\mathcal S}}   \def\CM{{\mathcal M}}   \def\CN{{\mathcal N}}
\def\CR{{\mathcal R}}   \def\CO{{\mathcal O}}   \def\CP{{\mathcal P}}
\def\CA{{\mathcal A}}   \def\CB{{\mathcal B}}   \def\CC{{\mathcal C}}
\def\CD{{\mathcal D}}   \def\CE{{\mathcal E}}   \def\CF{{\mathcal F}}
\def\CG{{\mathcal G}}   \def\CH{{\mathcal H}}   \def\CI{{\mathcal I}}
\def\CJ{{\mathcal J}}   \def\CK{{\mathcal K}}   \def\CL{{\mathcal L}}
\def\CT{{\mathcal T}}   \def\CU{{\mathcal U}}   \def\CV{{\mathcal V}}
\def\CZ{{\mathcal Z}}   \def\CX{{\mathcal X}}   \def\CY{{\mathcal Y}}
\def\CW{{\mathcal W}} \def\CQ{{\mathcal Q}} 
%%%%%
%%Macros for 'blackboard' letters (See (27) for display.)
\def\BBA {\mathbb A}   \def\BBb {\mathbb B}    \def\BBC {\mathbb C}
\def\BBD {\mathbb D}   \def\BBE {\mathbb E}    \def\BBF {\mathbb F}
\def\BBG {\mathbb G}   \def\BBH {\mathbb H}    \def\BBI {\mathbb I}
\def\BBJ {\mathbb J}   \def\BBK {\mathbb K}    \def\BBL {\mathbb L}
\def\BBM {\mathbb M}   \def\BBN {\mathbb N}    \def\BBO {\mathbb O}
\def\BBP {\mathbb P}   \def\BBR {\mathbb R}    \def\BBS {\mathbb S}
\def\BBT {\mathbb T}   \def\BBU {\mathbb U}    \def\BBV {\mathbb V}
\def\BBW {\mathbb W}   \def\BBX {\mathbb X}    \def\BBY {\mathbb Y}
\def\BBZ {\mathbb Z}

%%Macros for Ghotic (Fraktur) letters.
\def\GTA {\mathfrak A}   \def\GTB {\mathfrak B}    \def\GTC {\mathfrak C}
\def\GTD {\mathfrak D}   \def\GTE {\mathfrak E}    \def\GTF {\mathfrak F}
\def\GTG {\mathfrak G}   \def\GTH {\mathfrak H}    \def\GTI {\mathfrak I}
\def\GTJ {\mathfrak J}   \def\GTK {\mathfrak K}    \def\GTL {\mathfrak L}
\def\GTM {\mathfrak M}   \def\GTN {\mathfrak N}    \def\GTO {\mathfrak O}
\def\GTP {\mathfrak P}   \def\GTR {\mathfrak R}    \def\GTS {\mathfrak S}
\def\GTT {\mathfrak T}   \def\GTU {\mathfrak U}    \def\GTV {\mathfrak V}
\def\GTW {\mathfrak W}   \def\GTX {\mathfrak X}    \def\GTY {\mathfrak Y}
\def\GTZ {\mathfrak Z}   \def\GTQ {\mathfrak Q}

\font\Sym= msam10 % special symbols
\def\SYM#1{\hbox{\Sym #1}}
\newcommand{\bdw}{\prt\Gw\xspace}
\medskip
%%%%%%%%%%%%%%%%%%%%%%%%%%%%%%%%%%%%%%%%%%%%%%%%%%%%%%
\begin{abstract}
We study the boundary value problem with Radon measures  for nonnegative
solutions of $L_Vu:=-\Delta u+Vu=0$ in a bounded smooth domain
$\Gw$, when $V$ is a locally bounded nonnegative function. Introducing some specific capacity, we give sufficient conditions on a Radon measure $\gm$ on $\prt\Gw$ so that the problem can be solved. We study the reduced measure associated to this equation as well as the boundary trace of positive solutions.  In the appendix A. Ancona solves a question raised by M. Marcus and L. V\'eron concerning the vanishing set of the Poisson kernel of $L_V$ for an important class of potentials $V$. 
\end{abstract}
%\footnote{ \ 2010 MSC: , , 
\tableofcontents
%%%%%%%%%%%%%%%%%%%%%%%%%%%%%%%%%%%%%%%%%%%%%%%%%%%%%%%%%%%%%%%%%%%%%%%%%%%%%%%%%%%%%%%%%%%%%%%%%%%%%%%%%%%%%%%%%%%%%%%SECTION-THE REGULAR CASE%%%%%%%%%%%%%%%%%%%%%%%%%%%%%%%%%%%%%%%%%%%%%%%%%%%%%%%%%%%%%%%%%%%%%%%%%%%%%%%%%%%%%%%%%%%%%%%%%%%%%%%%%%%%%%%%%%%%%%%%%%%%%%%%%%%%%%
\mysection{Introduction}
 Let $\Gw$ be a smooth bounded domain of $\BBR^N$ and $V$ a locally bounded real valued measurable function defined in $\Gw$. The first question we adress is the solvability of the following non-homogeneous Dirichlet problem with a Radon measure for boundary data,
\begin{equation}\label{I1}
\left\{\BA{ll}
-\Gd u+Vu=0\qquad&\text{in }\Gw\\
\phantom{-\Gd u+V}
u=\gm\qquad&\text{in }\prt\Gw.
\EA\right.
\end{equation}
Let $\gf$ be the first (and positive) eigenfunction of $-\Gd$ in $W^{1,2}_0(\Gw)$. By a solution we mean a function $u\in L^1(\Gw)$, such that $Vu\in L^1_\gf$, which satisfies
 \begin{equation}\label{I2}
\myint{\Gw}{}\left(-u\Gd\gz+Vu\gz\right)dx=-\myint{\prt\Gw}{}\myfrac{\prt\gz}{\prt{\bf n}}d\gm.
\end{equation}
for any function $\gz\in C_0^1(\overline\Gw)$ such that $\Gd\gz\in L^\infty(\Gw)$. When $V$ is a bounded nonnegative function, it is straightforward that there exist a unique solution. However, it is less obvious to find general conditions which allow the solvability for any $\gm\in \frak M(\prt\Gw)$, the set of Radon measures on $\prt\Gw$. In order to avoid difficulties due to Fredholm type obstructions, we shall most often assume that $V$ is nonnegative, in which case there exists at most one solution. \smallskip

Let us denote by $K^\Gw$ the Poisson kernel in $\Gw$ and by $\BBK[\gm]$ the Poisson potential of a measure, that is
 \begin{equation}\label{I6}
\BBK[\gm](x):=\myint{\prt\Gw}{}K^\Gw(x,y)d\gm(y)\qquad\forall x\in\Gw.
\end{equation}
We first observe that, when $V\geq 0$ and the measure $\gm$ satisfies
 \begin{equation}\label{I5}
\myint{\Gw}{}\BBK[|\gm|](x)V(x)\gf(x)dx<\infty, 
\end{equation}
then problem $(\ref{I1})$ admits a solution. A Radon measure which satisfies $(\ref{I5})$ is called {\it an admissible measure} and a measure for which a solution exists is called {\it a good measure}. \smallskip

We first consider the {\it subcritical case} which means that the boundary value is solvable for any  $\gm\in \frak M(\prt\Gw)$. As a first result, we prove that any measure $\gm$ is admissible if $V$ is nonnegative and satisfies
 \begin{equation}\label{I4}
\sup_{\;\;\;\;\;\;y\in\prt\Gw}\!\!\!\!\! {\rm ess}\myint{\Gw}{}K^\Gw(x,y)V(x)\gf(x)dx<\infty,
\end{equation}
where $\gf$ is the first positive eigenfuntion of $-\Gd$ in $W^{1,2}_0(\Gw)$. Using estimates on the Poisson kernel, this condition is fulfilled if there exists $M>0$ such that for any $y\in\prt\Gw$,
 \begin{equation}\label{I3}
\myint{0}{D(\Gw)}\left(\myint{\Gw\cap B_r(y)}{}V(x)\gf^2(x) dx\right)\myfrac{dr}{r^{N+1}}\leq M
\end{equation}
where $D(\Gw)=diam (\Gw)$. We give also sufficient conditions which ensures that the boundary value problem $(\ref{I1})$ is stable from the weak*-topology of $\frak M(\prt\Gw)$ to $L^1(\Gw)\cap L^1_{V\gf}(\Gw)$. One of the sufficient conditions is that  $V\geq 0$ satisfies
\begin{equation}\label{stab0}
\lim_{\ge\to 0}\myint{0}{\ge}\left(\myint{\Gw\cap B_r(y)}{}V(x)\gf^2(x)dx\right)\myfrac{dr}{r^{N+1}}=0,
\end{equation} 
uniformly with respect to $y\in\prt\Gw$.
\smallskip

In the {\it supercritical case} problem $(\ref{I1})$ cannot be solved for any $\gm\in \frak M(\prt\Gw)$. In order to characterize positive good measures, we introduce a framework of nonlinear analysis which have been used by Dynkin and Kuznetsov (see \cite{Dbook1} and references therein) and Marcus and V\'eron \cite{MV1}  in their study of the  boundary value problems with measures
\begin{equation}\label{I8}
\left\{\BA{ll}
-\Gd u+|u|^{q-1}u=0\qquad&\text{in }\Gw\\
\phantom{-\Gd u+|u|^{q-1}}
u=\gm\qquad&\text{in }\prt\Gw,
\EA\right.
\end{equation}
where $q>1$. In these works, positive good measures on $\prt\Gw$ are completely characterized by the $C_{2/q,q'}$-Bessel in dimension N-1 and the following property:

\smallskip {\it A measure $\gm\in\frak M_+(\prt\Gw)$ is good for problem $(\ref{I8})$ if and only if it does charge Borel sets with zero $C_{2/q,q'}$-capacity, i.e
\begin{equation}\label{I9}
C_{2/q,q'}(E)=0\Longrightarrow\gm(E)=0\qquad\forall E\subset\prt\Gw,\,E\text{ Borel}.
\end{equation}
Moreover, any positive good measure is the limit of an increasing sequence $\{\gm_n\}$ of admissible measures which, in this case, are the positive measures belonging to the Besov space $B_{2/q,q'}(\prt\Gw)$.
They also characaterize  removable sets in terms of $C_{2/q,q'}$-capacity.}\medskip

In our present work, and always with $V\geq 0$, we use a capacity associated to the Poisson kernel $K^\Gw$ and which belongs to a class  studied by Fuglede \cite{Fu} \cite{Fu2}. It is defined by
\begin{equation}\label{I10}
C_V(E)=\sup\{\gm (E):\gm\in\frak M_+(\prt\Gw), \gm(E^c)=0,\,\norm{V\BBK[\gm]}_{L^1_\gf}\leq 1\},
\end{equation}
for any  Borel set $E\subset\prt\Gw$. Furtheremore $C_V(E)$ is equal to the value of its dual expression $C^*_V(E)$ defined by
\begin{equation}\label{I11}
C^*_V(E)=\inf\{\norm {f}_{L^\infty}:\check\BBK[f]\geq 1\quad\text {on }E\},
\end{equation}
where
\begin{equation}\label{I12}
\check\BBK[f](y)=\myint{\Gw}{}K^\Gw(x,y)f(x)V(x)\gf(x)dx\qquad\forall y\in\prt\Gw.
\end{equation}
If $E$ is a compact subset of $\prt\Gw$, this capacity is explicitely given by
\begin{equation}\label{I12'}
C_V(E)=C^*_V(E)=\max_{y\in E}\left(\myint{\Gw}{}K^\Gw(x,y)V(x)\gf(x) dx\right)^{-1}.
\end{equation}
\medskip

We denote by $Z_V$ the largest set with zero $C_V$ capacity, i.e.

\begin{equation}\label{I7}
Z_V=\left\{y\in \prt\Gw:\myint{\Gw}{}K^\Gw(x,y)V(x)\gf(x) dx=\infty\right\},
\end{equation}
and we prove the following.\smallskip 

\noindent{\it 1- If $\{\gm_n\}$ is an increasing sequence of positive good measures which converges to a measure $\gm$  in the weak* topology, then $\gm$ is a good measure. 
\smallskip 

\noindent 2- If $\gm\in \frak M_+(\prt\Gw)$ satisfies $\gm (Z_V)=0$, then $\gm$ is a good measure.\smallskip 

\noindent 3- A good measure $\gm$ vanishes on $Z_V$ if and only if there exists  an increasing sequence of positive admissible measures which converges to $\gm$  in the weak* topology.
}\medskip

In section 4 we study relaxation phenomenon in replacing $(\ref{I1})$ by the truncated problem
\begin{equation}\label{I1-k}
\left\{\BA{ll}
-\Gd u+V_ku=0\qquad&\text{in }\Gw\\
\phantom{-\Gd u+V_k}
u=\gm\qquad&\text{in }\prt\Gw.
\EA\right.
\end{equation}
where $\{V_{k}\}$ is an increasing sequence of positive bounded functions which converges to $V$
locally uniformly in $\Gw$. We adapt to the linear problem some of the principles of the reduced measure. This notion is introduced by Brezis, Marcus and Ponce \cite{BMP} in the study of the nonlinear Poisson equation
\begin{equation}\label{I10'}
-\Gd u+g(u)=\gm\qquad\text{in }\Gw\
\end{equation}
and extended to the Dirichlet problem 
\begin{equation}\label{I10''}
\left\{\BA{ll}
-\Gd u+g(u)=0\qquad&\text{in }\Gw\\
\phantom{-\Gd u+g()}
u=\gm\qquad&\text{in }\prt\Gw,
\EA\right.
\end{equation}
by Brezis and Ponce \cite{BP}.  In our construction, problem $(\ref{I1-k})$ admits a unique solution $u_k$. The sequence $\{u_k\}$ decreases and converges to some $u$ which satisfies a relaxed boundary value problem
\begin{equation}\label{I13}
\left\{\BA{ll}
-\Gd u+Vu=0\qquad&\text{in }\Gw\\
\phantom{-\Gd u+V}
u=\gm^*\qquad&\text{in }\prt\Gw.
\EA\right.
\end{equation}
The measure $\gm^*$ is called the {\it reduced measure} associated to $\gm$ and $V$. {\it Note that $\gm^*$ is the largest measure for which the problem
\begin{equation}\label{I15}
\left\{\BA{ll}
-\Gd u+Vu=0\qquad&\text{in }\Gw\\
\phantom{-\Gd u+V}
u=\gn\leq \gm\qquad&\text{in }\prt\Gw.
\EA\right.
\end{equation}
admits a solution}. This truncation process allows to construct the Poisson kernel $K^\Gw_V$ associated to the operator $-\Gd +V$ as being the limit of the decreasing limit of the sequence of kernel functions $\{K^\Gw_{V_k}\}$ asociated to $-\Gd +V_k$.
The solution $u=u_{\gm^*}$ of $(\ref{I13})$ is expressed by
\begin{equation}\label{I14}
u_{\gm^*}(x)=\myint{\prt\Gw}{}K^\Gw_V(x,y)d\gm(y)=\myint{\prt\Gw}{}K^\Gw_V(x,y)d\gm^*(y)\qquad\forall x\in\Gw.
\end{equation}
We define the vanishing set of $K^\Gw_V$ by
\begin{equation}\label{I16}
{\mathcal S}{\scriptstyle ing}_{_ V}(\Omega )=\{y\in \prt\Gw:K^\Gw_V(x_0,y)=0\},
\end{equation}
for some $x_0\in \Gw$, and thus for any $x\in \Gw$ by Harnack inequality. We prove 
\smallskip

\noindent 1- ${\mathcal S}{\scriptstyle ing}_{_ V}(\Omega)\subset Z_V$.
\smallskip

\noindent 2- $\gm^*=\gm\chi_{_{{\mathcal S}{\scriptstyle ing}_{_ V}(\Omega)}}$.\smallskip

\noindent {\it A challenging open problem is to give conditions on $V$ which imply ${\mathcal S}{\scriptstyle ing}_{_ V}(\Omega )= Z_V$}. \smallskip

The last section is devoted to the construction of the boundary trace of positive solutions of 
\begin{equation}\label{I19}
-\Gd u+Vu=0\qquad\text{in }\Gw,
\end{equation}
assuming $V\geq 0$. Using results of \cite{MV3}, we defined the regular set $\CR(u)$ of the boundary trace of $u$. This set is a relatively open subset of $\prt\Gw$ and the regular part of the boundary trace is represented by a positive Radon measure $\gm_u$ on $\CR(u)$. In order to study the singular set of the boundary trace $\CS(u):=\prt\Gw\setminus \CR(u)$, we adapt the sweeping method introduced by Marcus and V\'eron in \cite{MV4} for equation 
\begin{equation}\label{I20}
-\Gd u+g(u)=0\qquad\text{in }\Gw.
\end{equation}
If $\gm$ is a good positive measure concentrated on $\CS(u)$, and $u_\gm$ is the unique solution of $(\ref{I1})$ with boundary data $\gm$, we set $v_\gm=\min\{u,u_\gm\}$. Then $v_\gm$ is a positive super solution which admits a positive trace $\gg_u(\gm)\in\frak M_+(\prt\Gw)$. The extended boundary trace $Tr^e(u)$ of $u$ is defined by
\begin{equation}\label{I21}
\gn(u)(E):=Tr^e(u)(E)=\sup\{\gg_u(\gm)(E):\gm\text { good}, E\subset\prt\Gw,\,E\text{ Borel}\}.
\end{equation}
 {\it Then $Tr^e(u)$ is a Borel measure on $\Gw$. If we assume moreover that 
 \begin{equation}\label{I22}
\lim_{\ge\to 0}\myint{0}{\ge}\left(\myint{\Gw\cap B_r(y)}{}V(x)\gf^2(x)dx\right)\myfrac{dr}{r^{N+1}}=0\qquad\text{uniformly with respect to }y\in\prt\Gw,
\end{equation}
then $Tr^e(u)$ is a bounded measure and therefore a Radon measure. Finally, if $N=2$ and $(\ref{I22})$ holds, or if $N\geq 3$ and there holds
 \begin{equation}\label{I22'}
\lim_{\ge\to 0}\myint{0}{\ge}\left(\myint{\Gw\cap B_r(y)}{}V(x)(\gf(x)-\ge)^2_+dx\right)\myfrac{dr}{r^{N+1}}=0,
\end{equation}
uniformly with respect to $\ge\in (0,\ge_0] \text{ and }y\text{ s.t. }\gd_\Gw(x):=\dist(x,\prt\Gw)=\ge$, then $u=u_{\gn(u)}$.}\medskip

If $V(x)\leq v(\gf(x)$ for some $v$ which satisfies 
 \begin{equation}\label{I23}
\myint{0}{1}v(t)tdt<\infty,
\end{equation}
then Marcus and V\'eron proved in \cite{MV3} that $u=u_{\gn_u}$. Actually, when $V$ has such a geometric form, the assumptions $(\ref{I22})$-$(\ref{I22'})$ and $(\ref{I23})$ are equivalent. \smallskip

The Appendix, written by A. Ancona, answers a question raised by M. Marcus and L. V\'eron  in 2005 about the vanishing set of $K_V$ when $V$ is nonnegative and  $\gd_\Gw^2V$ is uniformly bounded. Such potentials play a very important role in the description of the fine trace of semilinear elliptic equations as in (\ref{I8}): actually, for such equations, $V=u^{q-1}$ satisfies this upper estimate as a consequence of Keller-Osserman estimate. The following result is proved\smallskip

\noindent{\it Let $y\in\prt\Gw$ and $C_{\ge,y}:=\{x\in\Gw:\gd_\Gw(x)\geq \ge|x-y|\}$ for $0<\ge<1$. If
 \begin{equation}\label{I24}
\myint{C_{\ge,y}}{}\myfrac{V(x) dx}{|x-y|^{N-2}}=\infty,
\end{equation}
for some $\ge>0$, then $y\in {\mathcal S}{\scriptstyle ing}_{_ V}(\Omega )$.}
%%%%%%%%%%%%%%%%%%%%%%%%%%%%%%%%%%%%%%%%%%%%%%%%%%%%%%%%%%%%%%%%%%%%%%%%%%%%%%%%%%%%%%%%%%%%%%%%%%%%%%%%%%%%%%%%%%%%%%%THE SUBCRITICAL CASE%%%%%%%%%%%%%%%%%%%%%%%%%%%%%%%%%%%%%%%%%%%%%%%%%%%%%%%%%%%%%%%%%%%%%%%%%%%%%%%%%%%%%%%%%%%%%%%%%%%%%%%%%%%%%%%%%%%%%%%%%%%%%%%%%%%%%%%%%%%%%%
\mysection{The subcritical case}
In the sequel $\Gw$ is a bounded smooth domain in $\BBR^N$ and $V\in L^\infty_{loc}$. We denote by $\gf$ the first eigenfunction of $-\Gd$ in $W^{1,2}_0(\Gw)$, $\gf>0$ with the corresponding eigenvalue $\gl$, by
$\frak M(\prt\Gw)$ the space of bounded Radon measures on $\prt\Gw$ and by $\frak M_+(\prt\Gw)$ its positive cone. For any positive Radon measure on $\prt\Gw$, we shall denote by the same symbol the corresponding outer regular bounded Borel measure. Conversely, for any outer regular bounded Borel $\gm$, we denote by the same expression $\gm$ the Radon measure defined on $C(\prt\Gw)$ by
$$\gz\mapsto \gm(\gz)=\myint{\prt\Gw}{}\gz d\gm.
$$
If $\gm\in \frak M(\prt\Gw)$, we are concerned with the following problem
\begin{equation}\label{bvp1}
\left\{\BA{ll}
-\Gd u+Vu=0\qquad&\text{in }\Gw\\
\phantom{-\Gd u+V}
u=\gm\qquad&\text{in }\prt\Gw.
\EA\right.
\end{equation}

\bdef{bvpdef} Let  $\gm\in \frak M(\prt\Gw)$. We say that $u$ is a weak solution of $(\ref{bvp1})$, if $u\in L^1(\Gw)$, $Vu\in L^1_\gf(\Gw)$ and, for any $\gz\in C^1_0(\overline\Gw)$ with $\Gd \gz\in L^\infty(\Gw)$, there holds
 \begin{equation}\label{bvp2}
\myint{\Gw}{}\left(-u\Gd\gz+Vu\gz\right)dx=-\myint{\prt\Gw}{}\myfrac{\prt\gz}{\prt {\bf n}}d\gm.
\end{equation}
\es
In the sequel we put
$$T(\Gw):=\{\gz\in C^1_0(\overline\Gw)\text { such that }\Gd \gz\in L^\infty(\Gw)\}.
$$

We recall the following estimates obtained by Brezis \cite{Br2}
\bprop{brez} Let  $\gm\in L^1(\prt\Gw)$ and $u$ be a weak solution of problem $(\ref{bvp1})$. Then there holds
 \begin{equation}\label{brez1}
\norm u_{L^1(\Gw)}+\norm {V_+u}_{L_\gf^1(\Gw)}\leq \norm {V_{-}u}_{L_\gf^1(\Gw)}+c\norm \gm_{L^1(\prt\Gw)}
\end{equation}
 \begin{equation}\label{brez2}
\myint{\Gw}{}\left(-|u|\Gd\gz+V|u|\gz\right)dx\leq-\myint{\prt\Gw}{}\myfrac{\prt\gz}{\prt {\bf n}}|\gm|dS
\end{equation}
and 
 \begin{equation}\label{brez3}
\myint{\Gw}{}\left(-u_+\Gd\gz+Vu_+\gz\right)dx\leq-\myint{\prt\Gw}{}\myfrac{\prt\gz}{\prt {\bf n}}\gm_+dS,
\end{equation}
for all $\gz\in T(\Gw)$, $\gz\geq 0$.
\es

We denote by $K^{\Gw}(x,y)$ the Poisson kernel in $\Gw$ and by $\BBK[\gm]$ the Poisson potential of $\gm\in \frak M(\prt\Gw)$ defined by
 \begin{equation}\label{Pois}
\BBK[\gm](x)=\myint{\prt\Gw}{}K^{\Gw}(x,y)d\gm(y)\qquad\forall x\in\Gw.
\end{equation}

\bdef {adm} A measure $\gm$ on $\prt\Gw$ is {\bf admissible} if
 \begin{equation}\label{adm}
\myint{\Gw}{}\BBK[|\gm|](x)|V(x)|\gf(x)dx<\infty.
\end{equation}
It is {\bf good} if problem $(\ref{bvp1})$ admits a weak solution.
\es

We notice that, if there exists at least one admissible positive measure $\gm$, then
\begin{equation}\label{glob}
\myint{\Gw}{}V(x)\gf^2(x)dx<\infty.
\end{equation}
%%%%%%%%%%%%%%%%%%%%%%%%%%%%%%%%%%%%%%%%%%%%%%%%%%%%%%%%%%%%%%%%%%%%%%%%%%%%%%%%%%%%%%%%%%%%%%%%%%%%%%%%%%%%%%
\bth {Exist} Assume $V\geq 0$, then problem $(\ref{bvp1})$ admits at most one solution.  Furthermore, if $\gm$ is admissible, then there exists a unique solution that we denote $u_\gm$.
\es
\Proof Uniqueness follows from $(\ref{brez1})$. For existence we can assume $\gm\geq 0$. For any $k\in\BBN_{*}$ set $V_{k}=\inf\{V,k\}$ and denote by $u:=u_{k}$ the solution of 
\begin{equation}\label{L4}\left\{\BA {ll}
-\Gd u+V_{k}(x)u=0\qquad&\text {in }\Gw\\
\phantom{-\Gd u+V(x)}
u=\gm\qquad&\text {on }\prt\Gw.
\EA\right.\end{equation} 
Then $0\leq u_{k}\leq \BBK[\gm]$. 
By the maximum principle, $u_{k}$ is decreasing and converges to some $u$, and 
$$0\leq V_{k}u_{k}\leq V\BBK[\gm]. 
$$
Thus, by dominated convergence theorem $V_{k}u_{k}\to Vu$ in $L^1_\gf$. Setting $\gz\in T(\Gw)$ and letting $k$ tend to infinity in  equality
\begin{equation}\label{L5}
\myint{\Gw}{}\left(-u_{k}\Gd\gz+V_{k}u_{k}\gz\right)dx=-\myint{\prt\Gw}{}\frac{\prt\gz}{\prt{\bf n}}d\gm, 
\end{equation}
implies that $u$ satisfies $(\ref{bvp2})$.\qeda\medskip

\noindent\Remark If $V$ changes sign, we can put $\tilde u=u+\BBK[\gm]$. Then $(\ref{bvp1})$ is equivalent to
\begin{equation}\label{bvp4}
\left\{\BA{ll}
-\Gd \tilde u+V\tilde u=V\BBK[\gm]&\qquad\text{in }\Gw\\
\phantom{-\Gd \tilde u+V}
\tilde u=0\qquad&\text{in }\prt\Gw.
\EA\right.
\end{equation}
This is a Fredholm type problem (at least if the operator $\phi\mapsto R(v):=(-\Gd)^{-1}(V\phi)$ is compact in $L^1_\gf(\Gw)$). Existence will be ensured by orthogonality conditions. \medskip

If we assume that $V\geq 0$ and
\begin{equation}\label{bvp5}
\myint{\Gw}{}K^\Gw(x,y)V(x)\gf(x)dx<\infty,
\end{equation}
for some $y\in \prt\Gw$, then $\gd_y$ is admissible. The following result yields to the  solvability of $(\ref{bvp1})$ for any $\gm\in\frak M_+(\Gw)$.
\bprop{Uncond} Assume $V\geq 0$ and the integrals $(\ref{bvp5})$ are bounded uniformly with respect to $y\in\prt\Gw$. Then any measure on $\prt\Gw$ is admissible.
\es
\Proof If $M$ is the upper bound of these integrals and $\gm\in \frak M_+(\prt\Gw)$, we have,
\begin{equation}\label{ener4}
\myint{\Gw}{}\BBK[\gm](x)V(x)\gf(x)dx=\myint{\prt\Gw}{}\left(\myint{\Gw}{}K^\Gw(x,y)V(x)\gf(x)dx\right)d\gm(y)\leq 
M\gm(\prt\Gw),
\end{equation}
 by Fubini's theorem. Thus $\gm$ is admissible.
\qeda\medskip

\noindent\Remark Since the Poisson kernel in $\Gw$ satisfies the two-sided estimate
\begin{equation}\label{ener7}
c^{-1}\myfrac{\gf(x)}{|x-y|^{N}}\leq K^\Gw(x,y)\leq c\myfrac{\gf(x)}{|x-y|^{N}}\qquad \forall (x,y)\in \Gw\ti\prt\Gw,
\end{equation}
for some $c>0$, assumption $(\ref{bvp5})$ is equivalent to 
\begin{equation}\label{ener6}
\myint{\Gw}{}\myfrac{V(x)\gf^2(x)}{|x-y|^{N}}dx<\infty.
\end{equation}
This implies $(\ref{glob})$ in particular. If we set $D_y=\max\{|x-y|:x\in\Gw\}$, then

$$\BA {l}\myint{\Gw}{}\myfrac{V(x)\gf^2(x)}{|x-y|^{N}}dx= 
\myint{0}{D_y}\left(\myint{\{x\in\Gw:|x-y|=r\}}{}\!\!\!\!\!\!V(x)\gf^2(x)dS_r(x)\right)\myfrac{dr}{r^{N}}\\[4mm]
\phantom{\myint{\Gw}{}\myfrac{V(x)\gf(x)}{|x-y|^{N-1}}dx}
=\displaystyle\lim_{\ge\to 0}\left(\left[r^{-N}\myint{\Gw\cap B_r(y)}{}\!\!\!\!\!\!V(x)\gf^2(x)dx\right]_\ge^{D_y}
+N\myint{\ge}{D_y}\left(\myint{\Gw\cap B_r(y)}{}\!\!\!\!\!\!V(x)\gf^2(x)dx\right)\myfrac{dr}{r^{N+1}}\right)
\EA$$
(both quantity may be infinite). Thus, if we assume 
 \begin{equation}\label{marc0}
\myint{0}{D_y}\left(\myint{\Gw\cap B_r(y)}{}\!\!\!\!\!\!V(x)\gf^2(x)dx\right)\myfrac{dr}{r^{N+1}}<\infty,
 \end{equation}
there holds 
 \begin{equation}\label{Marc0'}
 \liminf_{\ge\to 0}\ge^{-N}\myint{\Gw\cap B_\ge(y)}{}\!\!\!\!\!\!V(x)\gf^2(x)dS=0.
 \end{equation}
Consequently
 \begin{equation}\label{Marc0''}
\myint{\Gw}{}\myfrac{V(x)\gf^2(x)}{|x-y|^{N}}dx=
D_y^{-N}\myint{\Gw}{}V(x)\gf^2(x)dx+N\myint{0}{D_y}\left(\myint{\Gw\cap B_r(y)}{}\!\!\!\!\!\!V(x)\gf^2(x)dx\right)\myfrac{dr}{r^{N+1}}.
 \end{equation}
Therefore $(\ref{bvp5})$ holds and $\gd_y$ is admissible. \medskip

%%%%%%%%%%%%%%%%%%%%%%%%%%%%%%%%%%%%%%%%%%%%%%%%%%%%%%%%%%%%%%%%%%%STRENGTHEN%%%%%%%%%%%%%%%%%%%%%%%%%%%%%%%%%%%%%%%%%%%%%%%%%%%%%%%%%%%%%%%%%%%%%%%%%%%%%%%%%%%%%%%%%% 

As a natural extension of \rprop{Uncond}, we have the following stability result.
%%%%%%%%%%%%%%%%%%%%%%%%%%%%%%%%%%%%%%%%%%%%%%%%%%%%%%%%%%%%%%%%%%%%%%%%%%%%%%%%%%%%%%%%%%%%%%%%%%%%%%%%%%%%%%
\bth{stab} Assume $V\geq 0$ and 
\begin{equation}\label{cond'}
\lim_{\tiny\BA {l}E\text{ Borel}\\|E|\to 0\EA}\myint{E}{}K^\Gw(x,y)V(x)\gf(x)dx=0\quad\text {uniformly with respect to }y\in\prt\Gw.
\end{equation}
If $\gm_n$ is a sequence of positive Radon measures on $\prt\Gw$ converging to $\gm$ in the weak* topology, then $u_{\gm_n}$ converges to $u_{\gm}$ in $L^1(\Gw)\cap L^1_{V\gf} (\Gw)$ and locally uniformly in $\Gw$.
\es
\Proof 
We put $u_{\gm_n}:=u_n$. By the maximum principle $0\leq u_{n}\leq \BBK[\gm_n]$. Furthermore, it follows from $(\ref{brez1})$ that
 \begin{equation}\label{stab1}
\norm {u_{n}}_{L^1(\Gw)}+\norm {Vu_{n}}_{L_\gf^1(\Gw)}\leq c\norm {\gm_n}_{L^1(\prt\Gw)}\leq C.
\end{equation}
Since $-\Gd u_n$ is bounded in $L_\gf^1(\Gw)$, the sequence $\{u_n\}$ is relatively compact in $L^1(\Gw)$ by the regularity theory for elliptic equations. Therefore, there exist a subsequence  $u_{n_k}$ and some function $u\in L^1(\Gw)$ with $Vu\in L^1_\gf(\Gw)$ such that
$u_{n_k}$ converges to $u$ in $L^1(\Gw)$, almost everywhere on $\Gw$  and locally uniformly in $\Gw$ since $V\in L^\infty_{loc}(\Gw)$. The main question is to prove the convergence of $Vu_{n_k}$ in $L^1_\gf(\Gw)$. If $E\subset\Gw$ is any Borel set, there holds
$$\BA {l}
\myint{E}{}u_{n}V(x)\gf(x) dx\leq\myint{E}{}\BBK[\gm_n]V(x)\gf(x) dx\\[3mm]\phantom{\myint{E}{}u_{n}V(x)\gf(x) dx}
\leq\myint{\prt\Gw}{}\left(\myint{E}{}K^\Gw(x,y)V(x)\gf(x) dx\right)d\gm_n(y)\\[3mm]\phantom{\myint{E}{}u_{n}V(x)\gf(x) dx}
\leq M_n\displaystyle\max_{y\in\prt\Gw}\myint{E}{}K^\Gw(x,y)V(x)\gf(x) dx
,\EA$$
where $M_n:=\gm_n(\prt\Gw)$. Thus
 \begin{equation}\label{stab5}\BA {l}
\myint{E}{}u_{n}V(x)\gf(x) dx
\leq M_n\displaystyle\max_{y\in\prt\Gw}\myint{E}{}K^\Gw(x,y)V(x)\gf(x) dx.
\EA\end{equation}
Then, by $(\ref{cond'})$,
$$
\lim_{|E|\to 0}\myint{E}{}u_{n}V(x)\gf(x) dx=0.
$$
As a consequence the set of function $\{u_{n}\gf V\}$ is uniformly integrable. By Vitali's theorem
$Vu_{n_k}\to Vu$ in $L^1_\gf(\Gw)$. Since 
 \begin{equation}\label{stab6}
\myint{\Gw}{}\left(-u_n\Gd\gz+Vu_n\gz\right)dx=-\myint{\prt\Gw}{}\myfrac{\prt\gz}{\prt {\bf n}}d\gm_n,
\end{equation}
for any $\gz\in T(\Gw)$, the function $u$ satisfies $(\ref{bvp2})$.
\qeda\medskip

Assumption $(\ref{cond'})$ may be difficult to verify and the following result gives an easier formulation.

\bprop{CS} Assume $V\geq 0$ satisfies
 \begin{equation}\label{Marc0}
\lim_{\ge\to 0}\myint{0}{\ge}\left(\myint{\Gw\cap B_r(y)}{}V(x)\gf^2(x)dx\right)\myfrac{dr}{r^{N+1}}=0\quad\text{uniformly with respect to }y\in\prt\Gw.
\end{equation} 
Then $(\ref{cond'})$ holds.
\es
\Proof If $E\subset\Gw$ is a Borel set and $\gd>0$, we put $E_\gd=E\cap B_\gd(y)$ and $E^c_\gd=E\setminus E_\gd$. Then
$$ \myint{E}{}\myfrac{V(x)\gf^2(x)}{|x-y|^N}dx= \myint{E_\gd}{}\myfrac{V(x)\gf^2(x)}{|x-y|^N}dx+ \myint{E_\gd^c}{}\myfrac{V(x)\gf^2(x)}{|x-y|^N}dx.
$$
Clearly
 \begin{equation}\label{Marc1}
 \myint{E_\gd^c}{}\myfrac{V(x)\gf^2(x)}{|x-y|^N}dx\leq \gd^{-N} \myint{E.}{}V(x)\gf^2(x) dx.
\end{equation} 
Since  $(\ref{marc0})$ holds for any $y\in\prt\Gw$, $(\ref{Marc0''})$ implies
 \begin{equation}\label{Marc3}
\myint{E_\gd}{}\myfrac{V(x)\gf^2(x)}{|x-y|^N}dx=
\gd^{-N}\myint{E_\gd}{}V(x)\gf^2(x)dx+N\myint{0}{\gd}\left(\myint{E\cap B_r(y)}{}V(x)\gf^2(x)dx\right)\myfrac{dr}{r^{N+1}}.
\end{equation} 
Using $(\ref{Marc0})$, for any $\ge>0$, there exists $s_0>0$ such that for any  $s>0$ and $y\in\prt\Gw$
$$s\leq s_0\Longrightarrow N\myint{0}{s}\left(\myint{B_r(y)}{}V(x)\gf^2(x)dx\right)\myfrac{dr}{r^{N+1}}\leq\ge/2.
$$
We fix $\gd=s_0$. Since $(\ref{glob})$ holds,
\begin{equation}\label{cond+}
\lim_{\tiny\BA {l}E\text{ Borel}\\|E|\to 0\EA}\myint{E}{}V(x)\gf^2(x)dx=0.
\end{equation}
Then there exists $\eta>0$ such that for any Borel set $E\subset\Gw$,
$$|E|\leq\eta\Longrightarrow \myint{E}{}V(x)\gf^2(x)dx\leq s^N_0\ge/4.
$$
Thus
$$\myint{E}{}\myfrac{V(x)\gf^2(x)}{|x-y|^N}dx\leq\ge.
$$
This implies the claim by $(\ref{ener7})$.\qeda\medskip

An assumption which is used in \cite[Lemma 7.4]{MV3} in order to prove the existence of a boundary trace of any positive solution of $(\ref{I19})$ is that there exists some nonnegative measurable function $v$ defined on $\BBR_+$ such that
 \begin{equation}\label{t1}
\abs {V(x)}\leq v(\gf(x))\quad\forall x\in\Gw\quad\text{and }\myint{0}{s}tv(t)dt<\infty\quad\forall s>0.
\end{equation}

In the next result we show that condition $(\ref{t1})$ implies $(\ref{cond'})$.

\bprop{Tr} Assume $V$ satisfies $(\ref{t1})$. Then
\begin{equation}\label{condabs}
\lim_{\tiny\BA {l}E\text{ Borel}\\|E|\to 0\EA}\myint{E}{}K^\Gw(x,y)\abs{V(x)}\gf(x)dx=0\quad\text {uniformly with respect to }y\in\prt\Gw.
\end{equation}
\es
\Proof Since $\prt\Gw$ is $C^2$, there exist $\ge_0>0$ such that any for any $x\in\Gw$ satisfying $\gf(x)\leq\ge_0$, there exists a unique $\gs(x)\in \prt\Gw$ such that $|x-\gs(x)|=\gf(x)$. 
We use $(\ref{Marc0})$ in \rprop{CS} under the equivalent form
 \begin{equation}\label{Marc0+}
\lim_{\ge\to 0}\myint{0}{\ge}\left(\myint{\Gw\cap C_r(y)}{}|V(x)|\gf^2(x)dx\right)\myfrac{dr}{r^{N+1}}=0\quad\text{uniformly with respect to }y\in\prt\Gw,
\end{equation} 
in which we have replaced $ B_r(y)$ by the the cylinder $C_r(y):=\{x\in\Gw:\gf(x)<r,|\gs(x)-y|<r\}$. Then
$$\BA {l}
\myint{0}{\ge}\left(\myint{\Gw\cap C_r(y)}{}|V(x)|\gf^2(x)dx\right)\myfrac{dr}{r^{N+1}}
\leq c\myint{0}{\ge}\left(\myint{0}{r}v(t)t^2dt\right)\myfrac{dr}{r^{2}}\\[4mm]
\phantom{\myint{0}{\ge}\left(\myint{\Gw\cap C_r(y)}{}|V(x)|\gf^2(x)dx\right)\myfrac{dr}{r^{N+1}}}
\leq  c\myint{0}{\ge}v(t)\left(1-\myfrac{t}{\ge}\right)tdt\\[4mm]
\phantom{\myint{0}{\ge}\left(\myint{\Gw\cap C_r(y)}{}|V(x)|\gf^2(x)dx\right)\myfrac{dr}{r^{N+1}}}
\leq c\myint{0}{\ge}v(t)tdt.
\EA$$
Thus $(\ref{Marc0})$ holds.\qeda\medskip

 %%%%%%%%%%%%%%%%%%%%%%%%%%%%%%%%%%%%%%%%%%%%%%%%%%%%%%%%%%%%%%%%%%%%%%%%%%%%%%%%%%%%%%%%%%%%%%%%%%%%%%%%%%%%%%%%%%%%%%%SECTION-THE CAPACITARY APPROACH%%%%%%%%%%%%%%%%%%%%%%%%%%%%%%%%%%%%%%%%%%%%%%%%%%%%%%%%%%%%%%%%%%%%%%%%%%%%%%%%%%%%%%%%%%%%%%%%%%%%%%%%%%%%%%%%%%%%%%%%%%%%%%%%%%%%%%

\mysection{The capacitary approach}

Throughout this section $V$ is a locally bounded nonnegative and measurable function defined on $\Gw$. We assume that there exists a positive measure $\gm_0$ on $\prt\Gw$ such that 
\begin{equation}\label{capa1}
\myint{\Gw}{}\BBK[\gm_0]V(x)\gf(x)dx=\CE(1,\gm_0)<\infty.
\end{equation}

\bdef {energy}  If $\gm\in\frak M_+(\prt\Gw)$ and $f$ is a nonnegative measurable function defined in $\Gw$ such that 
$$(x,y)\mapsto \BBK[\gm](y)f(x)V(x)\gf(x)\in L^1(\Gw\ti\prt\Gw;dx\otimes d\gm),$$
we set
\begin{equation}\label{capa2}
\CE(f,\gm)=\myint{\Gw}{}\left(\myint{\prt\Gw}{}K^\Gw(x,y)d\gm(y)\right)f(x)V(x)\gf(x)dx.
\end{equation}
\es
If we put
\begin{equation}\label{capa3}
\check\BBK_V[f](y)=\myint{\Gw}{}K^\Gw(x,y)f(x)V(x)\gf(x)dx,
\end{equation}
then, by Fubini's theorem, $\check\BBK_{V}[f]<\infty$, $\gm$-almost everywhere on $\prt\Gw$ and
\begin{equation}\label{capa4}
\CE(f,\gm)=\myint{\prt\Gw}{}\left(\myint{\Gw}{}K^\Gw(x,y)f(x)V(x)\gf(x)dx\right)d\gm(y).
\end{equation}
\bprop {lsc} Let $f$ be fixed. Then \smallskip

\noindent (a) $y\mapsto \check \BBK_V[f](y)$ is lower semicontinuous on $\prt\Gw$.\smallskip

\noindent (b) $\gm\mapsto \CE(f,\gm)$ is lower semicontinuous on $\frak M_+(\prt\Gw)$ in the weak*-topology
\es
\Proof Since $y\mapsto K^\Gw(x,y)$ is continuous, statement (a) follows by Fatou's lemma. If ${\gm_n}$ is a sequence in $\frak M_+(\prt\Gw)$ converging to some $\gm$ in the weak*-topology, then 
$\BBK[\gm_n]$ converges to $\BBK[\gm]$ everywhere in $\Gw$. By Fatou's lemma 
$$\CE(f,\gm)\leq \liminf_{n\to\infty}\myint{\Gw}{}\BBK[\gm_n](x)f(x)V(x)\gf(x)dx=\liminf_{n\to\infty}\CE(f,\gm_n).
$$
\qeda

Notice that if $V\gf f\in L^p(\Gw)$, for $p>N$, then $\BBG[Vf\gf]\in C^1(\overline\Gw)$ and
 \begin{equation}\label{capa5}
\check \BBK[f](y):=\myint{\Gw}{}K^\Gw(x,y)V(x)f(x)\gf(x)dx=-\myfrac{\prt}{\prt{\bf n}}\BBG[Vf\gf](y).
\end{equation}
This is in particular the case if $f$ has compact support in $\Gw$.\medskip

\bdef{mspace} We denote by $\frak M^V(\prt\Gw)$ the set of all measures $\gm$ on $\prt\Gw$ such that
$V\BBK[\gm]\in L^1_{\gf}(\Gw)$. If $\gm$ is such a measure, we denote
 \begin{equation}\label{mspace1}
\norm\gm_{\frak M^V}= \myint{\Gw}{}\abs{\BBK[\gm](x)}V(x)\gf(x) dx=\norm{V\BBK[\gm]}_{L^1_\gf}.
\end{equation}
\es

Clearly $\norm{\,.\,}_{\frak M^V}$ is a norm. The space  $\frak M^V(\prt\Gw)$ is not complete but its positive cone $\frak M_+^V(\prt\Gw)$ is complete. If $E\subset\prt\Gw$ is a Borel subset, we put
$$\frak M_+(E)=\{\gm\in\frak M_+(\prt\Gw): \gm (E^c)=0\}\quad\text{and }\;\frak M^V_+(E)=\frak M_+(E)\cap \frak M^V(\prt\Gw).
$$
\bdef{cap} If $E\subset\prt\Gw$ is any Borel subset we set
 \begin{equation}\label{capa6}
C_V(E):=\sup \{\gm (E) :\gm\in\frak M^V_+(E), \norm\gm_{\frak M^V}\leq 1\}.
\end{equation}
\es
We notice that $(\ref{capa6})$ is equivalent to
 \begin{equation}\label{fug1}
C_V(E):=\sup \left\{\myfrac{\gm (E)}{\norm\gm_{\frak M^V}} :\gm\in\frak M^V_+(E)\right\}.
\end{equation}

\bprop{capaci} The set function $C_V$ satisfies. 
 \begin{equation}\label{fu0}
C_V(E)\leq\sup_{y\in E} \left(\myint{\Gw}{}K^\Gw(x,y)V(x)\gf(x)dx\right)^{-1}\quad\forall E\subset\prt\Gw,\, E\text{ Borel},
\end{equation}
and equality holds in $(\ref{fu0})$ if $E$ is compact. Moreover, 
 \begin{equation}\label{fug0}
C_V(E_1\cup E_2)=\sup \{C_V(E_1),C_V(E_2)\}\quad\forall E_i\subset\prt\Gw,\, E_i\text{ Borel}.
\end{equation}
\es
\Proof Notice that $E\mapsto C_V(E)$ is a nondecreasing set function for the inclusion relation and that $(\ref{capa6})$ implies
 \begin{equation}\label{fu1}
\gm (E)\leq C_V(E)\norm\gm_{\frak M^V}\qquad\forall\gm\in\frak M^V_+(E).
\end{equation}
Let $E\subset\prt\Gw$ be a Borel set and $\gm\in \frak M_+(E)$. Then
$$\BA {l}\norm {\gm}_{\frak M^V}=\myint{E}{}\left(\myint{\Gw}{}K^{\Gw}(x,y)V(x)\gf(x) dx\right)d\gm(y)\\[4mm]
\phantom{\norm {\gm_{\frak M^V}}}\geq \gm(E)\displaystyle{\inf_{y\in E}\myint{\Gw}{}K^{\Gw}(x,y)V(x)\gf(x) dx}.
\EA$$
Using $(\ref{capa6})$ we derive
 \begin{equation}\label{fu2}
C_V(E)\leq \sup_{y\in E}\left(\myint{\Gw}{}K^{\Gw}(x,y)V(x)\gf(x) dx\right)^{-1}.
\end{equation}
If $E$ is compact, there exists $y_0\in E$ such that
$$\inf_{y\in E}\myint{\Gw}{}K^{\Gw}(x,y)V(x)\gf(x) dx=\myint{\Gw}{}K^{\Gw}(x,y_0)V(x)\gf(x) dx,
$$
since $y\mapsto \check \BBK[1](y)$ is l.s.c.. Thus
$$\norm {\gd_{y_0}}_{\frak M^V}=\gd_{y_0}(E)\myint{\Gw}{}K^{\Gw}(x,y_0)V(x)\gf(x) dx
$$
and 
$$C_V(E)\geq \myfrac{\gd_{y_0}(E)}{\norm {\gd_{y_0}}_{\frak M^V}}=\sup_{y\in E}\left(\myint{\Gw}{}K^{\Gw}(x,y)V(x)\gf(x) dx\right)^{-1}.
$$
Therefore equality holds in $(\ref{fu0})$.
%%%%%%%%%%%%%%%%%%%%%%%%%%%%%%%%%%%%%%%%%%%%%%%%%%%%%%%%%%%%%%%%%%%%%%%%%%%%%%%%%%%%%%%%%%%%%%%%%%%%%%%%%%%%
Identity $(\ref{fug0})$ follows  $(\ref{fu0})$ when there is equality. Moreover it holds if $E_1$ and $E_2$ are two arbitrary compact sets. Since $C_V$ is eventually an inner regular capacity (i.e. $C_V(E)=\sup \{C_V(K): K\subset E,\, K\text { compact}\}$) it holds for any Borel set. However we give below a self-contained proof.
If  $E_1$ and $E_2$ be two disjoint Borel subsets of $\prt\Gw$, for any $\ge>0$ 
there exists $\gm\in \frak M^V_+(E_1\cup E_2)$ such that 
$$\myfrac{\gm (E_1)+\gm (E_2)}{\norm{\gm}_{\frak M^V}}\leq C_V(E_{1}\cup E_{2})\leq \myfrac{\gm (E_1)+\gm (E_2)}{\norm{\gm}_{\frak M^V}}+\ge.
$$
Set $\gm_i=\chi_{_{E_i}}\gm$. Then $\gm_i\in \frak M^V_+(E_i)$ and  $\norm{\gm}_{\frak M^V}=\norm{\gm_1}_{\frak M^V}+\norm{\gm_2}_{\frak M^V}$. By $(\ref{fu1})$ 
 \begin{equation}\label{capa6*}
 C_V(E_{1}\cup E_{2})\leq \myfrac{\norm{\gm_1}_{\frak M^V}}{\norm{\gm_1}_{\frak M^V}+\norm{\gm_2}_{\frak M^V}}C_V(E_{1})+\myfrac{\norm{\gm_2}_{\frak M^V}}{\norm{\gm_1}_{\frak M^V}+\norm{\gm_2}_{\frak M^V}}C_V(E_{2})+\ge
\end{equation}
This implies that there exists $\gth\in [0,1]$ such that
 \begin{equation}\label{capa6**}
 C_V(E_{1}\cup E_{2})\leq \gth C_V(E_{1})+(1-\gth)C_V(E_{2})\leq \max\{C_V(E_{1}),C_V(E_{2})\}.
\end{equation}
Since $C_V(E_{1}\cup E_{2})\geq\max\{C_V(E_{1}),C_V(E_{2})\}$ as $C_V$ is increasing,
 \begin{equation}\label{capa6***}
  E_1\cap E_2=\emptyset\Longrightarrow C_V(E_{1}\cup E_{2})=\max\{C_V(E_{1}),C_V(E_{2})\}.
\end{equation}
If  $E_1\cap E_2\neq\emptyset$, then $E_1\cup E_2=E_1\cup (E_2\cap E^c_1)$
and therefore
$$C_V(E_{1}\cup E_{2})=\max\{C_V(E_{1}),C_V(E_2\cap E^c_1)\}\leq \max\{C_V(E_{1}),C_V(E_{2})\}.
$$
Using again $(\ref{fug1})$ we derive $(\ref{fug0})$.
\qeda\medskip

The following set function is the dual expression of $C_V(E)$.

\bdef {dual}For any Borel set $E\subset\prt\Gw$, we set
 \begin{equation}\label{capa7}
C^*_V(E):=\inf \{\norm f_{L^\infty} :\check \BBK[f](y)\geq 1\quad\forall y\in E\}.
\end{equation}
\es 
The next result is stated in \cite[p 922]{Fu2} using minimax theorem and the fact that $K^\Gw$ is lower semi continuous in $\Gw\ti\prt\Gw$. Although the proof is not explicited, a simple adaptation of the proof of \cite[Th 2.5.1]{AH} leads to the result. 

\bprop{Equi} For any compact set  $E\subset\prt\Gw$,
\begin{equation}\label {zeroset}
C_V(E)= C^*_V(E).
 \end{equation}
 \es

In the same paper \cite{Fu2}, formula $(\ref{fu0})$ with equality is claimed (if $E$ is compact).

\bth{G} If $\{\gm_{n}\}$ is an increasing sequence of good measures converging to some measure $\gm$ in the weak* topology, then $\gm$ is good.
\es
\Proof We use formulation $(\ref{D9})$. We take for test function the function $\eta$ solution of
\begin{equation}\left\{\BA {ll}
-\Gd \eta=1\qquad&\text{in }\Gw\\
\phantom{-\Gd}\eta=0\qquad&\text{on }\Gw,
\EA\right.\end{equation}
there holds
$$\myint{\Gw}{}\left(1+V\right)u_{\gm_{n}}\eta dx=
-\myint{\prt\Gw}{}\myfrac{\prt\eta}{\prt{\bf n}}d\gm_{n}\leq c^{-1}\gm_{n}(\prt\Gw)\leq c^{-1}\gm(\prt\Gw)
$$
where $c>0$ is such that
$$c^{-1}\geq -\myfrac{\prt\eta}{\prt{\bf n}}\geq c\quad\text {on }\prt\Gw.
$$
Since $\{u_{\gm_{n}}\}$ is increasing and $\eta\leq c\gf$ by Hopf boundary lemma, we can let $n\to\infty$ by the monotone convergence theorem. If $u:=\lim_{n\to\infty}u_{\gm_{n}}$, we obtain
$$\myint{\Gw}{}\left(1+V\right)u \eta dx\leq c^{-1}\gm(\prt\Gw).
$$
Thus $u$ and $\gf Vu$ are in $L^1(\Gw)$. Next, if 
$\gz\in C_{0}^{1}(\overline\Gw)\cap C^{1,1}(\overline\Gw)$, then
$u_{\gm_{n}}|\Gd\gz|\leq Cu_{\gm_{n}}$ and $Vu_{\gm_{n}}|\gz|\leq CVu_{\gm_{n}}\eta$. Because the sequence $\{u_{\gm_{n}}\}$ and 
$\{Vu_{\gm_{n}}\eta\}$ are uniformly integrable, the same holds for 
$\{u_{\gm_{n}}\Gd\gz\}$ and $\{Vu_{\gm_{n}}\gz\}$. Considering
$$
\myint{\Gw}{}\left(-u_{\gm_{n}}\Gd\gz+Vu_{\gm_{n}}\gz\right)dx=-\myint{\prt\Gw}{}\myfrac{\prt\gz}{\prt{\bf n}}d\gm_{n}.
$$
it follows by Vitali's theorem, 
$$
\myint{\Gw}{}\left(-u\Gd\gz+Vu\gz\right)dx=-\myint{\prt\Gw}{}\myfrac{\prt\gz}{\prt{\bf n}}d\gm.
$$
Thus $\gm$ is a good measure.\qeda\medskip

We define the {\it singular boundary set} $Z_V$ by
\begin{equation}\label{Z0}
Z_V=\left\{y\in\prt\Gw:\myint{\Gw}{}K^\Gw(x,y)V(x)\gf(x) dx=\infty\right\}.
\end{equation}
Since $\check\BBK[1]$ is l.s.c., it is a Borel function and $Z_V$ is a Borel set. The next result characterizes  the good measures.
\bprop{van} Let $\gm$ be an admissible positive measure. Then $\gm(Z_V)=0$.
\es
\Proof If $K\subset Z_V$ is compact,  $\gm_K=\chi_{_K}\gm$ is admissible, thus, by Fubini theorem
$$\norm{\gm_K}_{\frak M^V}=\myint{K}{}\left(\myint{\Gw}{}K^{\Gw}(x,y)V(x)\gf(x)dx\right)d\gm(y)<\infty.
$$
Since 
$$\myint{\Gw}{}K^{\Gw}(x,y)V(x)\gf(x)dx\equiv \infty\qquad\forall y\in K
$$
it follows that $\gm(K)=0$. This implies $\gm(Z_V)=0$ by regularity.\qeda
%%%%%%%%%%%%%%%%%%%%%%%%%%%%%%%%%%%%%%%%%%%%%%THEOREM%%%%%%%%%%%%%%%%%%%%%%%%%%%%%%%%%%%%%%%%%%%%%%%%%%%%%%%%%
\bth{approx} Let $\gm\in \frak M_+(\prt\Gw)$ such that
\begin{equation}\label{Z1}
\gm(Z_V)=0.
\end{equation}
Then $\gm $ is good.
\es
\Proof Since $\check\BBK[1]$ is l.s.c., for any $n\in\BBN_*$, 
$$K_n:=\{y\in \prt\Gw:\check\BBK[1](y)\leq n\}
$$
is a compact subset of $\prt\Gw$. Furthermore $K_n\cap Z_V=\emptyset$ and $\cup K_n=Z_V^c$. Let 
$\gm_n=\chi_{_{K_n}}\gm$, then
\begin{equation}\label{Z4}
\CE(1,\gm_n)=\myint{\Gw}{}\BBK[\gm_n]V(x)\gf(x)dx\leq n\gm_n(K_n).
\end{equation}
Therefore $\gm_n$ is admissible. By the monotone convergence theorem, $\gm_n\uparrow \chi_{_{Z_{V^c}}}\gm$ and by \rth{G}, $\chi_{_{Z_{V^c}}}\gm$ is good. Since $(\ref{Z1})$ holds, $\chi_{_{Z_{V^c}}}\gm=\gm$, which ends the proof.
\qeda\medskip

The full characterization of the good measures in the general case appears to be difficult without any further assumptions on $V$. However the following  holds
\bth{equiv} Let $\gm\in \frak M_+(\prt\Gw)$ be a good measure. The following assertions are equivalent:\smallskip

\noindent (i) $\gm(Z_V)=0$.\smallskip

\noindent (ii) There exists an increasing sequence of admissible measures $\{\gm_n\}$ which converges to $\gm$ in the weak*-topology.
\es
%%%%%%%%%
\Proof If (i) holds, it follows from the proof of \rth{approx} that the sequence $\{\gm_n\}$ increases and converges to $\gm$. If (ii) holds, any admissible measure $\gm_n$ vanishes on $Z_V$ by \rprop{van}. Since
$\gm_n\leq \gm$, there exists an increasing sequence of $\gm$-integrable functions $h_n$ such that $\gm_n=h_n \gm$. Then
$\gm_n(Z_V)$ increases to $\gm(Z_V)$ by the monotone convergence theorem. The conclusion follows from the fact that $\gm_n(Z_V)=0$. \qeda

%%%%%%%%%%%%%%%%%%%%%%%%%%%%%%%%%%%%%%%%%%%%%%%%%%%%%%%%%%%%%%%%%%%%%%%%%%%%%%%%%%%%%%%%%%%%%%%%%%%%%%%%%%%%%%%%%%%%%%%%%%%SECTION%%REPRESENTATION%%FORMULA%%%%%%%%%%%%%%%%%%%%%%%%%%%%%%%%%%%%%%%%%%%%%%%%%%%%%%%%%%%%%%%%%%%%%%%%%%%%%%%%%%%%%%%%%%%%%%%%%%%%%%%%%%%%%%%%%%%%%%%%%%%%%%%%%%%%%

\mysection{Representation formula and reduced measures}

We recall the construction of the Poisson kernel for $-\Gd +V$: if we look for a solution of 
 \begin{equation}\label{D1}\left\{
 \BA {ll}
 -\Gd v+V(x)v=0\qquad&\text {in }\Gw\\ 
 \phantom{ -\Gd v+V(x)}
 v=\gn\qquad&\text {in }\prt\Gw,
 \EA\right.
\end{equation}
where $\gn\in\frak M(\prt\Gw)$, $V\geq 0$, $V\in L^\infty_{loc}(\Gw)$, we can consider an increasing sequence of smooth domains 
$\Gw_n$ such that $\overline\Gw_n\subset\Gw_{n+1}$ and $\cup_n\Gw_n=\cup_n\overline\Gw_n=\Gw$. For each of these domains, denote by $K^{\Gw}_{V\chi_{_{\Gw_n}}}$ the Poisson kernel of $-\Gd+V\chi_{_{\Gw_n}}$ in $\Gw$ and by $\BBK_{V\chi_{_{\Gw_n}}}[.]$ the corresponding operator. We denote by $K^\Gw:=K^{\Gw}_0$ the Poisson kernel in $\Gw$ and by $\BBK[.]$ the Poisson operator in $\Gw$. Then the solution $v:=v_n$ of
 \begin{equation}\label{D2}\left\{
 \BA {ll}
 -\Gd v+V\chi_{_{\Gw_n}}v=0\qquad&\text {in }\Gw\\ 
 \phantom{  -\Gd v+V\chi_{_{\Gw_n}}}
 v=\gn&\text {in }\prt\Gw,
 \EA\right.
\end{equation}
is expressed by
 \begin{equation}\label{D3}
 v_n(x)=\myint{\prt\Gw}{}K^{\Gw}_{V\chi_{_{\Gw_n}}}(x,y)d\gn (y)=\BBK_{V\chi_{_{\Gw_n}}}[\gn](x).
\end{equation}
If $G^\Gw$ is the Green kernel of $-\Gd$ in $\Gw$ and $\BBG[.]$ the corresponding Green operator, $(\ref{D3})$ is equivalent to 
 \begin{equation}\label{D4}
v_n(x)+\myint{\Gw}{}G^\Gw(x,y)(V\chi_{_{\Gw_n}}v_n)(y)dy=\myint{\prt\Gw}{}K^{\Gw}(x,y)d\gn (y),
\end{equation}
equivalently
$$v_n+\BBG[V\chi_{_{\Gw_n}}v_n]=\BBK[\gn].
$$
Notice that this equality is equivalent to the weak formulation of problem $(\ref{D2})$: for any $\gz\in T(\Gw)$, there holds
 \begin{equation}\label{D4'}
\myint{\Gw}{}\left(-v_n\Gd\gz+V\chi_{_{\Gw_n}}v_n\gz\right)dx=-\myint{\prt\Gw}{}\myfrac{\prt\gz}{\prt{\bf n}}d\gn.
\end{equation}
Since $n\mapsto K^{\Gw}_{V\chi_{_{\Gw_n}}}$ is decreasing, the sequence $\{v_n\}$ inherits this property and there exists 
 \begin{equation}\label{D5}
\lim_{n\to\infty}K^{\Gw}_{V\chi_{_{\Gw_n}}}(x,y)=K^{\Gw}_{V}(x,y).
\end{equation}
By the monotone convergence theorem, 
 \begin{equation}\label{D6}  
\lim_{n\to\infty}v_n(x)=v(x)=\myint{\prt\Gw}{}K^{\Gw}_{V}(x,y)d\gn(y).
 \end{equation}
By Fatou's theorem
 \begin{equation}\label{D7}
 \myint{\Gw}{}G^\Gw(x,y)V(y)v(y)dy\leq \liminf_{n\to\infty}\myint{\Gw}{}G^\Gw(x,y)(V\chi_{_{\Gw_n}}v_n)(y)dy,
\end{equation}
and thus,  
 \begin{equation}\label{D8}
v(x)+ \myint{\Gw}{}G^\Gw(x,y)V(y)v(y)dy\leq\BBK[\gn](x)\qquad\forall x\in \Gw.
\end{equation}
Now the main question is to know whether $v$ keeps the boundary value $\gn$. Equivalently, whether the equality holds in $(\ref{D7})$ with $\lim$ instead of $\liminf$, and therefore in $(\ref{D8})$. This question is associated to the notion of reduced measured in the sense of Brezis-Marcus-Ponce: since $Vv\in L^1_{\gf}(\Gw)$ and 
 \begin{equation}\label{D9}
 -\Gd v+V(x)v=0\qquad\text {in }\Gw
 \end{equation}
 holds, the function $v+\BBG[Vv]$ is positive and harmonic in $\Gw$. Thus it admits a boundary trace $\gn^*\in \frak M_{+}(\prt\Gw)$ and
 \begin{equation}\label{D10}
v+ \BBG[Vv]=\BBK[\gn^*].
\end{equation}
Equivalently $v$  satisfies the relaxed problem
 \begin{equation}\label{D11}\left\{
 \BA {ll}
 -\Gd v+V(x)v=0&\qquad\text {in }\Gw\\ 
 \phantom{ -\Gd v+V(x)}
 v=\gn^*\qquad&\text {in }\prt\Gw,
 \EA\right.
 \end{equation}
 and thus $v=u_{\gn^*}$. Noticed that $\gn^*\leq\gn$ and the mapping $\gn\mapsto\gn^*$ is nondecreasing.
\bdef{red}The measure $\gn^*$ is the {\it reduced measure} associated to $\gn$.\es
%%%%%%%%%%%%%%%%%%%%%%%%%%%%%%%%%%%%%%%%%%%%%%%%%%%%%%%%%%%%%%%%%%%%%%%%%%%%%%%%%%%%%%%%%%%%%%%%%%%%%%%%%%%%%%%
\bprop {opt} There holds $\BBK_V[\gn]=\BBK_V[\gn^*]$. Furthermore
the reduced measure $\gn^*$ is the largest measure  for which the following problem  
\begin{equation}\label{D12}\left\{
 \BA {ll}
  \phantom{.} -\Gd v+V(x)v=0\qquad&\text {in }\Gw\\ 
  \phantom{}
\!\gl\in \frak M_{+}(\prt\Gw),\;\gl\leq\gn
 \\ 
 \phantom{ -\Gd v+V(x);.}
 v=\gl\qquad&\text {in }\prt\Gw,
 \EA\right.
 \end{equation}
 admits a solution.
\es
\Proof The first assertion follows from the fact that $v=\BBK_V[\gn]$ by $(\ref{D5})$ and $v=u_{\gn^*}=\BBK_V[\gn^*]$ by $(\ref{D11})$.
It is clear that $\gn^*\leq\gn$ and that the problem $(\ref{D12})$ admits a solution for $\gl=\gn^*$. If $\gl$ is a positive measure smaller than 
$\gm$, then $\gl^*\leq \gm^*$. But if there exist some $\gl$ such that the problem $(\ref{D12})$ admits a solution, then $\gl=\gl^*$. This implies the claim.\qeda\medskip

As a consequence of the characterization of $\gn^*$ there holds

\bcor{Uni} Assume $V\geq 0$ and let $\{V_k\}$ be an increasing sequence of nonnegative bounded measurable functions converging to $V$ a.e. in $\Gw$. Then the solution $u_k$ of 
 \begin{equation}\label{uni1}\left\{\BA {ll}
-\Gd u+V_ku=0\qquad&\text{in }\Gw\\ \phantom{-\Gd u+V_k}
u=\gn\qquad&\text{in }\prt\Gw,
\EA\right.\end{equation}
converges to $u_{\gn^*}$.
\es
\Proof The previous construction shows that 
$u_k=\BBK_{V_k}[\gn]$ decreases to some $\tilde u$ which satisfies a relaxed equation, the boundary data of which, $\tilde\gn^*$,  is the largest measure $\gl\leq\gn$ for which problem $(\ref{D12})$ admits a solution. Therefore $\tilde\gn^*=\gn^*$ and $\tilde u=u_{\gn^*}$. Similarly  $\{K^\Gw_{V_k}\}$  decreases and converges to $K^\Gw_{V}$. \qeda\medskip

We define the {\it boundary vanishing set of $K^\Gw_V$}  by

 \begin{equation}\label{van1}
{\mathcal S}{\scriptstyle ing}_{_ V}(\Omega ):=\{y\in\prt\Gw\,|\,K^\Gw_V(x,y)=0\}\quad\text{for some } x\in\Gw.
\end{equation}
Since $V\in L^\infty_{loc}(\Gw)$, ${\mathcal S}{\scriptstyle ing}_{_ V}(\Omega )$ is independent of $x$ by Harnack inequality; furthermore it is a Borel set. This set is called the set of {\it finely irregular boundary points} by E. B. Dynkin; the reason for such a denomination will appear in the Appendix.

\bth{vanish}  Let $\gn\in\frak M_+(\prt\Gw)$.\smallskip

\noindent (i) If $\gn(({\mathcal S}{\scriptstyle ing}_{_ V}(\Omega ))^c)=0$, then $\gn^*=0$.\smallskip

\noindent (ii) There always holds ${\mathcal S}{\scriptstyle ing}_{_ V}(\Omega )\subset Z_V$. 
\es
\Proof The first assertion is clear since $\gn=\chi_{_{{\mathcal S}{\scriptstyle ing}_{_ V}(\Omega )}}\gn+\chi_{_{{\mathcal S}{\scriptstyle ing}_{_ V}(\Omega ))^c}}\gn=\chi_{_{{\mathcal S}{\scriptstyle ing}_{_ V}(\Omega )}}\gn$ and, by \rprop{opt},
$$u_{\gn^*}(x)=\BBK_V[\gn^*](x)=\myint{{\mathcal S}{\scriptstyle ing}_{_ V}(\Omega )}{}K^\Gw_V(x,y)d\gn(y)=0\qquad\forall x\in\Gw,
$$
by definition of ${\mathcal S}{\scriptstyle ing}_{_ V}(\Omega )$. For proving (ii), we assume that $C_V({\mathcal S}{\scriptstyle ing}_{_ V}(\Omega ))>0$; there exists $\gm\in\frak M^V_+( {\mathcal S}{\scriptstyle ing}_{_ V}(\Omega ))$ such that $\gm({\mathcal S}{\scriptstyle ing}_{_ V}(\Omega ))>0$. 
Since $\gm$ is admissible let $u_\gm$ be the solution of $(\ref{I1})$. Then $\gm^*=\gm$, thus $u_\gm=\BBK^V[\gm]$
and
$$\BBK^V[\gm](x)=\myint{\prt\Gw}{}K^\Gw_V(x,y)d\gm(y)=\myint{{\mathcal S}{\scriptstyle ing}_{_ V}(\Omega )}{}K^\Gw_V(x,y)d\gm(y)=0,
$$
contradiction. Thus $C_V({\mathcal S}{\scriptstyle ing}_{_ V}(\Omega ))=0$. Since $(\ref{fu0})$  implies that $Z_V$ is the largest Borel set with zero $C_V$-capacity, it implies ${\mathcal S}{\scriptstyle ing}_{_ V}(\Omega )\subset Z _V$.\qeda

 %%%%%%%%%%%%%%%%%%%%%%%%%%%%%%%%%%%%%%%%%%%%%%%%%%%%%%%%%%%%%%%%%%%%%%%%%%%%%%%%%%%%%%%%%%%%%%%%%%%%%%%%%%%%%%

%%%%%%%%%%%ANCONA%%%%%%%%%%%%%%%%%%%%%%%%%%%%%%

\medskip

In order to obtain more precise informations on ${\mathcal S}{\scriptstyle ing}_{_ V}(\Omega )$ some minimal regularity assumptions on $V$ are needed. We also recall the following result due to Ancona \cite{anc4} and developed in the appendix of the present work.
\bth{An} Assume $V\geq 0$ satisfies $\gd_\Gw^2V\in L^\infty(\Gw)$. If for some $y\in \prt\Gw$ and some cone 
$C_{y}$ with vertex $y$ such that 
$\overline {C}_{y}\cap B_r(y)\subset \Gw\cup\{y\}$ for some $r>0$
there holds
 \begin{equation}\label{D15}
 \myint{C_{y}}{}\myfrac{V(x)}{|x-y|^{N-2}}dx=\infty,
\end{equation}
then
 \begin{equation}\label{D16}
K_V^{\Gw}(x,y)=0\qquad\forall x\in \Gw.
\end{equation}
\es
This means that (\ref{D15}) implies that $y$ belongs to ${\mathcal S}{\scriptstyle ing}_{_ V}(\Omega )$. Set $\gd_\Gw(x)=\dist (x,\prt\Gw)$. We define the {\it conical  singular boundary set}
 \begin{equation}\label{D17}
\tilde Z_V=\left\{y\in \prt\Gw:\myint{C{\ge,y}}{}K^\Gw(x,y)V(x)\gf(x) dx=\infty\;\text{for some } \ge>0\right\}
\end{equation}
where $C_{\ge,y}:=\{x\in\Gw:\gd_\Gw(x)\geq \ge|x-y|\}$. Clearly $\tilde Z_V\subset  Z_V$.  
\bcor{Osc} Assume $V\geq 0$ satisfies $\gd_\Gw^2V\in L^\infty(\Gw)$.
Then $\tilde Z_V\subset{\mathcal S}{\scriptstyle ing}_{_ V}(\Omega )$.
\es
\Proof Let $y\in \tilde Z_V$. Since there exists $c>0$ such that 
 \begin{equation}\label{D18}
 c^{-1}V(x)|x-y|^{2-N}\leq K^\Gw(x,y)V(x)\gf(x)\leq cV(x)|x-y|^{2-N}\qquad\forall x\in C_{\ge,y}
 \end{equation}
the result follows immediately from (\ref{D15}), (\ref{D17}).\qeda \medskip

\noindent \Remark In situations coming from the nonlinear equation $-\Gd u+|u|^{q-1}u=0$ in $\Gw$ with $q>1$, $V=|u|^{q-1}$ not only satisfies $gd_\Gw^2V\in L^\infty(\Gw)$ but also the restricted oscillation condition: for any $y\in\prt\Gw$ and any open cone $C_y$ with vertex $y$ such that 
$C_y\Subset\Gw$, there exists $c>0$ such that
 \begin{equation}\label{D19}
 \forall (x,z)\in C_y\ti C_y, |x-y|=|z-y|\Longrightarrow c^{-1}\leq \myfrac{V(x)}{V(z)}\leq c.
\end{equation}
It is a consequence of the Keller-Osserman estimate and Harnack inequality. In this case condition 
(\ref{D15}) is equivalent to
 \begin{equation}\label{D20}
 \myint{0}{1}V(\gamma (t))tdt=\infty,
\end{equation}
at least for one path $\gamma\in C^{0,1}([0,1])$ such that $\gamma (0)=y$ and $\gamma((0,1]\subset C_y$ for some cone $C_y\Subset\Gw$.

%%%%%%%%%%%%%%%%%%%%%%%%%%%%%%%%%%%%%%%%%%%%%%%%%%%%%%%%%%%%%%%%%%%%%%CAPACITY%%%%%%%%%%%%%%%%%%%%%%%%%%%%%%%%%%%%%%%%%%%%%%%%%%%%%%%%%%%%%%%%%

%%%%%%%%%%%%%%%%%%%%%%%%%%%%%%%%%%%%%%%%%%%%%%%%%%%%%%%%%%%%%%%%%%%%%%%%
%%%%%%%%%%%%%%%%%%%%%%%%%%%%%%%%%%%%%%%%%%%%%%%%%%%%%%%%%%%%%%%%%%%%%%%%%%%%%%%%%%%%%%%%%%%%%%%%%%%%%%%%%%%%%%%%%%%%%%%%%%%SECTION%%BOUNDARY%%TRACE%%%%%%%%%%%%%%%%%%%%%%%%%%%%%%%%%%%%%%%%%%%%%%%%%%%%%%%%%%%%%%%%%%%%%%%%%%%%%%%%%%%%%%%%%%%%%%%%%%%%%%%%%%%%%%%%%%%%%%%%%%%%%%%%%%%%%
\mysection{The boundary trace}
%%%%%%%%%%%%%%%%%%%%%%%%%%%%%%%%%%%%%%%%%%%%%%%%%%%%%%%%%%%%%%%%%%%%%%%%%%%%%%%%%%%%%%%%%%%%%%%%%%%%%%%%%%%%%%%%%%%%%%%%%%%%SUB-SECTION%%%%%%%%%%%%%%%%%%%%%%%%%%%%%%%%%%%%%%%%%%%%%%%%%%%%%%%%%%%%%%%%%%%%%%%%%%%%%%%%%%%%%%

\subsection{The regular part}
In this section, $V\in L^\infty_{loc}(\Gw)$ is nonnegative. If $0<\ge\leq \ge_0$, we denote $\gd_{\Gw}(x)=\dist (x,\prt\Gw)$ for $x\in\Gw$, and set $\Gw_\ge:=\{x\in\Gw:\gd_{\Gw}(x)>\ge\}$, $\Gw'_\ge=\Gw\setminus\Gw_\ge$ and $\Gs_\ge=\prt\Gw_\ge$. It is well known that there exists $\ge_0$ such that, for any $0<\ge\leq \ge_0$  and any $x\in\Gw'_\ge$ there exists a unique projection $\gs(x)$ of $x$ on $\prt\Gw$ and any $x\in \Gw'_\ge$ can be written in a unique way under the form
$$x=\gs(x)-\gd_{\Gw}(x){\bf n}
$$
where $\bf n$ is the outward normal unit vector to $\prt\Gw$ at $\gs(x)$.
The mapping $x\mapsto(\gd_{\Gw}(x),\gs(x))$ is a $C^2$ diffeomorphism from $\Gw'_\ge$ to $(0,\ge_0]\ti\prt\Gw$. We recall the following definition given in \cite{MV3}. If $\CA$ is a Borel subset of $\prt\Gw$, we set 
$\CA_\ge=\{x\in \Gs_\ge:\gs(x)\in A\}$.
\bdef{Trdef} Let $\CA$ be a relatively open subset of $\prt\Gw$, $\{\gm_\ge\}$ be a set of Radon measures on $\CA_{\ge}$ $(0<\ge\leq \ge_0)$ and $\gm\in\frak M(\CA)$. We say that $\gm_\ge\rightharpoonup\gm$ in the weak*-topology if, for any $\gz\in C_c(\CA)$, 
 \begin{equation}\label{Y1}
\lim_{\ge\to 0}\myint{\CA_{\ge}}{}\gz(\gs(x))d\gm_\ge(x)=\myint{\CA}{}\gz d\gm.
\end{equation}
A function $u\in C(\Gw)$ possesses a boundary trace $\gm\in \frak M(\CA)$ if
 \begin{equation}\label{Y2}
\lim_{\ge\to 0}\myint{\CA_{\ge}}{}\gz(\gs(x))u(x)dS(x)=\myint{\CA}{}\gz d\gm\qq\forall \gz\in C_c(\CA).
\end{equation}
\es
The following result is proved in \cite[p 694]{MV3}.
\bprop{Trlem1} Let $u\in C(\Gw)$ be a positive solution of 
 \begin{equation}\label{F1}
-\Gd u+V(x)u=0\qquad\text{in }\Gw.
\end{equation}
Assume that, for some $z\in\prt\Gw$, there exists an open neighborhood $U$ of $z$ such that 
 \begin{equation}\label{F2}
\myint{U\cap\Gw}{}Vu\gf(x)dx<\infty.
\end{equation}
Then $u\in L^1(K\cap\Gw)$ for any compact subset $K\subset G$ and there exists a positive Radon measure $\gm$ on $\CA=U\cap\prt\Gw$ such that
 \begin{equation}\label{Y3}
\lim_{\ge\to 0}\myint{U\cap\Gs_\ge}{}\gz(\gs(x))u(x)dS(x)=\myint{\CA}{}\gz d\gm\qq\forall \gz\in C_c(U\cap\Gw).
\end{equation}
\es

Notice that any continuous solution of $(\ref{F1})$ in $\Gw$ belongs to  $W^{2,p}_{loc}(\Gw)$ for any $(1\leq p<\infty)$. This previous result yields to a natural definition of the regular boundary points.

\bdef {Trdefreg}Let $u\in C(\Gw)$ be a positive solution of $(\ref{F1})$. 
A point $z\in\prt\Gw$ is called a regular boundary point for $u$ if there exists an open neighborhood $U$ of $z$ such that $(\ref{Y3})$ holds. The set of regular boundary points is a relatively open subset of $\prt\Gw$, denoted by $\CR(u)$. The set $\CS(u)=\prt\Gw\setminus\CR(u)$  is the singular boundary set of $u$. It is a closed set. 
\es

By \rprop{Trlem1} and using a partition of unity, we see that there exists a positive Radon measure $\gm:=\gm_u$
on $\CR(u)$ such that $(\ref{Y3})$ holds with $U$ replaced by $\CR(u)$. The couple $(\gm_u,\CS(u))$ is called the {\bf boundary trace of $u$}.  {\it The main question of the boundary trace problem is to analyse the behaviour of $u$ near the set $\CS(u)$.} \medskip

For any positive good measure $\gm$ on $\prt\Gw$, we denote by $u_\gm$ the solution of $(\ref{D1})$ defined by $(\ref{D9})$-$(\ref{D10})$. 

\bprop{Trlem2} Let  $u\in C(\Gw)\cap W^{2,p}_{loc}(\Gw)$ for any $(1\leq p<\infty)$ be a positive solution of  
$(\ref{F1})$ in $\Gw$ with boundary trace $(\gm_u,\CS(u))$. Then $u\geq u_{\gm_u}$.
\es
%%%%
\Proof Let $G\subset \prt\Gw$ be a relatively open subset such that $\overline G\subset \CR(u)$
 with a $C^2$ relative boundary $\prt^* G=\overline G\setminus G$. There exists an increasing sequence of 
$C^2$ domains $\Gw_n$ such that $\overline G\subset\prt\Gw_n$, $\prt\Gw_n\setminus \overline G\subset\Gw$ and $\cup_n\Gw_n=\Gw$. For any $n$, let $v:=v_n$ be the solution of 
 \begin{equation}\label{Z1}\left\{\BA {ll}
-\Gd v+Vv=0\qq&\text{in }\Gw_n\\ \phantom{-\Gd v+V}
v=\chi_{_G}\gm\qq&\text{in }\prt\Gw_n.
\EA\right.\end{equation}
Let $u_n$ be the restriction of $u$ to $\Gw_n$. Since $u\in C(\Gw)$ and $Vu\gf\in L^1(\Gw_n)$, there also holds
$Vu\gf_n\in L^1(\Gw_n)$ where we have denoted by $\gf_n$ the first eigenfunction of $-\Gd$ in $W^{1,2}_0(\Gw_n)$. Consequently $u_n$ admits a regular boundary trace $\gm_n$ on $\prt\Gw_n$ (i.e. $\CR(u_n)=\prt\Gw_n$) and $u_n$ is the solution of 
 \begin{equation}\label{Z1'}\left\{\BA {ll}
-\Gd v+Vv=0\qq&\text{in }\Gw_n\\ \phantom{-\Gd v+V}
v=\gm_n\qq&\text{in }\prt\Gw_n.
\EA\right.\end{equation}
Furthermore $\gm_n|_G=\chi_{_G}\gm_{u}$. It follows from Brezis estimates and in particular $(\ref{brez3})$ that 
$u_n\leq u$ in $\Gw_n$. Since  $\Gw_n\subset\Gw_{n+1}$, $v_n\leq v_{n+1}$. Moreover 
$$v_n+\BBG^{\Gw_n}[Vv_n]=\BBK^{\Gw_n}[\chi_{_G}\gm]\qquad\text{in }\Gw_n.
$$
Since $\BBK^{\Gw_n}[\chi_{_G}\gm_{u}]\to \BBK^{\Gw}[\chi_{_G}\gm_{u}]$, and the Green kernels $G^{\Gw_n}(x,y)$ are increasing with $n$, it follows from monotone convergence that $v_n\uparrow v$ and there holds
$$v+\BBG^{\Gw}[Vv]=\BBK^{\Gw}[\chi_{_G}\gm_{u}]\qquad\text{in }\Gw.
$$
Thus $v=u_{\chi_{_G}\gm_{u}}$ and $u_{\chi_{_G}\gm_{u}}\leq u$. We can now replace $G$ by a sequence $\{G_k\}$ of relatively open sets with the same properties as $G$, $\overline G_k\subset G_k$ and $\cup_k G_k=\CR(u)$. Then $\{u_{\chi_{_{G_k}}\gm_{u}}\}$ is increasing and converges to some $\tilde u$. Since
$$u_{\chi_{_G{_k}}\gm_{u}}+\BBG^{\Gw}[Vu_{\chi_{_G{_k}}\gm_{u}}]=\BBK^{\Gw}[\chi_{_G{_k}}\gm_{u}],
$$
and $\BBK^{\Gw}[\chi_{_G{_k}}\gm]\uparrow \BBK^{\Gw}[\gm_{u}]$, we derive
$$\tilde u+\BBG^{\Gw}[V\tilde u]=\BBK^{\Gw}[\gm_{u}].
$$
This implies that $\tilde u=u_{\gm_{u} }\leq u$.\qeda
%%%%%%%%%%%%%%%%%%%%%%%%%%%%%%%%%%%%%%%%%%%%%%%%%%%%%%%%%%%%%%%%%%%%%%%%%%%%%%%%%%%%%%%%%%%%%%%%%%%%%%%%%%%%%%%%%%%%%%%%%%%%SUB-SECTION%%%%%%%%%%%%%%%%%%%%%%%%%%%%%%%%%%%%%%%%%%%%%%%%%%%%%%%%%%%%%%%%%%%%%%%%%%%%%%%%%%%%%%

\subsection{The singular part}
The following result is essentially proved in \cite[Lemma 2.8]{MV3}.

\bprop{Trlem3} Let $u\in C(\Gw)$  for any $(1\leq p<\infty)$ be a positive solution of $(\ref{F1})$ and suppose that $z\in \CS(u)$ and that there exists an open neighborhood $U_0$ of $z$ such that 
$u\in L^1(\Gw\cap U_0)$. Then for any open neighborhood $U$ of $z$, there holds
 \begin{equation}\label{Y4}
\lim_{\ge\to 0}\myint{U\cap\Gs_\ge}{}\gz(\gs(x))u(x)dS(x)=\infty.
\end{equation}
\es

As immediate consequences, we have

\bcor {tot1}Assume $u$ satisfies the regularity assumption of \rprop{Trlem2}.  Then for any $z\in\CS(u)$ and any open neighborhood $U$ of $z$, there holds
 \begin{equation}\label{Y5}
\limsup_{\ge\to 0}\myint{U\cap\Gs_\ge}{}\gz(\gs(x))u(x)dS(x)=\infty.
\end{equation}
\es

\bcor {tot2}Assume $u$ satisfies the regularity assumption of \rprop{Trlem2}. If $u\in L^1(\Gw)$,  Then for any $z\in\CS(u)$ and any open neighborhood $U$ of $z$, $(\ref{Y4})$ holds.
\es

The two next results give conditions on $V$ which imply that $\CS(u)=\emptyset$.

\bth{regN=2} Assume $N=2$, $V$ is nonnegative and satisfies $(\ref{cond'})$. If $u$ is a positive solution of  $(\ref{F1})$, then  $\CR(u)=\prt\Gw$.
\es
\Proof We assume that 
\bel{R1}
\myint{\Gw}{}V\gf udx=\infty.
\ee
If $0<\ge\leq \ge_0$, we denote by $(\gf_\ge,\gl_\ge)$ are the normalized first eigenfunction and first eigenvalue of $-\Gd$ in $W^{1,2}_0(\Gw_\ge)$, then
\bel{R2}\lim_{\ge\to 0}\myint{\Gw_\ge}{}V\gf_\ge udx=\infty.
\ee
Because
$$\myint{\Gw_\ge}{}(\gl_\ge+\gf_\ge V)u dx=-\myint{\prt\Gw_\ge}{}\myfrac{\prt\gf_\ge}{\prt{\bf n}}udS,
$$
and
$$c^{-1}\leq-\myfrac{\prt\gf_\ge}{\prt{\bf n}}\leq c,
$$
for some $c>1$ independent of $\ge$, there holds
\bel{R3}\lim_{\ge\to 0}\myint{\prt\Gw_\ge}{}udS=\infty.
\ee
Denote by $m_\ge$ this last integral and set $v_\ge=m^{-1}_\ge u$ and $\gm_\ge=m^{-1}_\ge u|_{\prt\Gw_\ge}$.
Then
\bel{R3'}v_\ge+\BBG^{\Gw_\ge}[Vv_\ge]=\BBK^{\Gw_\ge}[\gm_\ge]\qquad\text{in }\Gw_\ge
\ee
where 
\bel{R4}\BBK^{\Gw_\ge}[\gm_\ge](x)=\myint{\prt\Gw_\ge}{}K^{\Gw_\ge}(x,y) \gm_\ge(y)dS(y)\ee
 is the Poisson potential of $\gm_\ge$ in $\Gw_\ge$ and
$$\BBG^{\Gw_\ge}[Vu](x)=\myint{\Gw_\ge}{}G^{\Gw_\ge}(x,y)V(y)u(y)dy,
$$
the Green potential of $Vu$ in $\Gw_\ge$. Furthermore
\bel{R5}\left\{\BA {ll}
-\Gd v_\ge+Vv_\ge=0\qquad&\text{in }\Gw_\ge\\\phantom{-\Gd v_\ge+V}
v_\ge=\gm_\ge\qquad&\text{in }\prt\Gw_\ge.
\EA\right.\ee
By Brezis estimates and regularity theory for elliptic equations, $\{\chi_{_{\Gw_\ge}}v_\ge\}$ is relatively compact in $L^1(\Gw)$ and in the local uniform topology of $\Gw_{\ge}$. Up to a subsequence $\{\ge_n\}$, $\gm_{\ge_n}$ converges to a probability measure $\gm$ on $\prt\Gw$ in the weak*-topology. It is classical that
$$\BBK^{\Gw_{\ge_n}}[\gm_{\ge_n}]\to \BBK[\gm]
$$
locally uniformly in $\Gw$, and $\chi_{_{\Gw_{\ge_n}}}v_{\ge_n}\to v$  in the local uniform topology of  $\Gw$,  and a.e. in $\Gw$. Because $G^{\Gw_\ge}(x,y)\uparrow G^{\Gw}(x,y)$, there holds for any $x\in \Gw$
\bel{R6}\lim_{n\to\infty}\chi_{_{\Gw_{\ge_n}}}(y)G^{\Gw_{\ge_n}}(x,y)V(y)v_{\ge_n}(y)=G^{\Gw}(x,y)V(y)v(y) \quad\text{for almost all  }y\in \Gw
\ee
Furthermore $v_{\ge_n}\leq \BBK^{\Gw_{\ge_n}}[\gm_{\ge_n}]$ reads
$$v_{\ge_n}(y)\leq c\gf_{\ge_n}(y)
\myint {\prt\Gw_n}{}\myfrac{\gm_{\ge_n}(z)dS(z)}{|y-z|^2}.
$$ 
In order to go to the limit in the expression
\bel{R7}L_n:=\BBG^{\Gw_{\ge_n}}[Vv_{\ge_n}](x)=\myint{\Gw}{}\chi_{_{\Gw_{\ge_n}}}(y)G^{\Gw_{\ge_n}}(x,y)V(y)v_{\ge_n}(y)dy,
\ee
we may assume that $x\in \Gw_{\ge_1}$ where $0<\ge_1\leq\ge_0$ is fixed and write 
$\Gw=\Gw_{\ge_1}\cup\Gw'_{\ge_1}$ where 
$$\Gw'_{\ge_1}=\Gw\setminus\Gw_{\ge_1}:=\{x\in\Gw:\dist (x,\prt\Gw)\leq\ge_1\}$$
and $L_n=M_n+P_n$ where
\bel{R8}M_n=\myint{\Gw_{\ge_1}}{}\chi_{_{\Gw_{\ge_n}}}(y)G^{\Gw_{\ge_n}}(x,y)V(y)v_{\ge_n}(y)dy
\ee
and
\bel{R9}P_n=\myint{\Gw'_{\ge_1}}{}\chi_{_{\Gw_{\ge_n}}}(y)G^{\Gw_{\ge_n}}(x,y)V(y)v_{\ge_n}(y)dy.
\ee
Since
$$\BA {l}
\chi_{_{\Gw_{\ge_1}}}(y)G^{\Gw_{\ge_n}}(x,y)V(y)v_{\ge_n}(y)
\leq c\chi_{_{\Gw_{\ge_1}}}(y)\abs{\ln(|x-y|)}V(y)v_{\ge_n}(y)\\[2mm]\phantom{\chi_{_{\Gw_{\ge_1}}}(y)G^{\Gw_{\ge_n}}(x,y)V(y)v_{\ge_n}(y)}
\leq c\norm{V}_{L^\infty(\Gw_{\ge_1})}\chi_{_{\Gw_{\ge_1}}}(y)\abs{\ln(|x-y|)}v_{\ge_n}(y),
\EA$$
it follows by the dominated convergence theorem that
\bel{R10}\lim_{n\to\infty}M_n=\myint{\Gw_{\ge_1}}{}G^{\Gw}(x,y)V(y)v(y)dy.
\ee
Let $E\subset\Gw$ be a Borel subset. Then $G^{\Gw_{\ge_n}}(x,y)\leq c(x)\gf_{\ge_n}(y)$ if $y\in \Gw'_{\ge_1}$. By Fubini,
\bel{R11}\BA {l}
\myint{\Gw'_{\ge_1}\cap E}{}\chi_{_{\Gw_{\ge_n}}}(y)G^{\Gw_{\ge_n}}(x,y)V(y)v_{\ge_n}(y)dy
\leq cc(x)\myint {\prt\Gw_n}{}\!\!\!\left(\myint{\Gw'_{\ge_1}\cap E}{}\chi_{_{\Gw_{\ge_n}}}\!\!(y)\myfrac{\gf^2_{\ge_n}(y)V(y)}{|y-z|^2}dy\right) \gm_{\ge_n}(z)dS(z)\\[4mm]
\phantom{\myint{A_{\ge_1}\cap E}{}\chi_{_{\Gw_{\ge_n}}}(y)G^{\Gw_{\ge_n}}(x,y)V(y)v_{\ge_n}(y)dy}
\leq cc(x)\displaystyle\max_{z\in\prt\Gw_{\ge_n}}\myint{\Gw'_{\ge_1}\cap E}{}\chi_{_{\Gw_{\ge_n}}}(y)\myfrac{\gf^2_{\ge_n}(y)V(y)}{|y-z|^2}dy
\EA\ee
If $y\in \Gw_{\ge_n}\cap E$, there holds $\gf(y)=\gf_{\ge_n}(y)+\ge_n$. If $z\in  \prt\Gw_{\ge_n}\cap E$ and we denote by $\gs(z)$ the projection of $z$ onto $\prt\Gw$, there holds $|y-\gs(z)|\leq |y-z|+\ge_n$. By monotonicity
\begin{equation}\label{ineq}
\myfrac{\gf_{\ge_n}(y)}{|y-z|}\leq \myfrac{\gf_{\ge_n}(y)+\ge_n}{|y-z|+\ge_n}\leq 
\myfrac{\gf(y)}{|y-\gs(z)|},
\end{equation}
thus
\begin{equation}\label{ineq2}\BA {l}
\myint{\Gw'_{\ge_1}\cap E}{}\chi_{_{\Gw_{\ge_n}}}(y)G^{\Gw_{\ge_n}}(x,y)V(y)v_{\ge_n}(y)dy
\leq cc(x)\displaystyle\max_{z\in\prt\Gw}\myint{\Gw'_{\ge_1}\cap E}{}\chi_{_{\Gw_{\ge_n}}}(y)\myfrac{\gf^2(y)V(y)}{|y-z|^2}dy.
\EA\end{equation}
By $(\ref{cond'})$ this last integral goes to zero if $\abs{\Gw'_{\ge_1}\cap E\cap \Gw_{\ge_n}}\to 0 $. Thus by Vitali's theorem, the sequence of functions $\{\chi_{_{\Gw_{\ge_n}}}(.)G^{\Gw_{\ge_n}}(x,.)V(y)v_{\ge_n}(.)\}_{n\in\BBN}$ is uniformly integrable in $y$, for any $x\in\Gw$. It implies that
\begin{equation}\label{ineq3}
\lim_{n\to\infty}\myint{\Gw}{}\chi_{_{\Gw_{\ge_n}}}(y)G^{\Gw_{\ge_n}}(x,y)V(y)v_{\ge_n}(y)dy
=\myint{\Gw}{}G^\Gw(x,y)V(y)v(y)dy,
\end{equation}
and there holds $v+\BBG[Vv]=\BBK[\gm]$. Since $u=m_\ge v_\ge$ in $\Gw$ and $m_\ge\to\infty$, we get a contradiction since it would imply $u\equiv\infty$. 
\qeda\medskip

In order to deal with the case $N\geq 3$ we introduce an additionnal assumption of stability.

\bth{regN>2} Assume $N\geq 3$. Let $V\in L^\infty_{loc}(\Gw)$, $V\geq 0$ such that
\bel{3-1}
\lim_{\tiny\BA {l}E\text{ Borel}\\|E|\to 0\EA}\myint{E}{}V(y)\myfrac{(\gf(y)-\ge)_+^2}{|y-z|^{N}}dy=0\quad\text {uniformly with respect to }z\in\Gs_\ge \text { and }\ge\in (0,\ge_0].
\ee
 If $u$ is a positive solution of  $(\ref{F1})$, then  $\CR(u)=\prt\Gw$.
\es
\Proof We proceed as in \rth{regN=2}. All the relations $(\ref{R1})$-$(\ref{R10})$ are valid and $(\ref{R11})$ has to be replaced by 
\bel{3-2}\BA {l}
\myint{\Gw'_{\ge_1}\cap E}{}\chi_{_{\Gw_{\ge_n}}}(y)G^{\Gw_{\ge_n}}(x,y)V(y)v_{\ge_n}(y)dy
\leq cc(x)\displaystyle\max_{z\in\Gs_{\ge_n}}\myint{\Gw'_{\ge_1}\cap E}{}\chi_{_{\Gw_{\ge_n}}}(y)\myfrac{\gf^2_{\ge_n}(y)V(y)}{|y-z|^{N+1}}dy.
\EA\ee
Since $(\ref{ineq})$ is no longer valid, $(\ref{ineq})$ is replaced by
\begin{equation}\label{3-3}\BA {l}
\myint{\Gw'_{\ge_1}\cap E}{}\chi_{_{\Gw_{\ge_n}}}(y)G^{\Gw_{\ge_n}}(x,y)V(y)v_{\ge_n}(y)dy
\leq cc(x)\displaystyle\max_{z\in\Gs_{\ge_n}}\myint{ E}{}V(y)\myfrac{(\gf(y)-\ge_n)_+^2}{|y-z|^{N+1}}dy.
\EA\end{equation}
By $(\ref{3-1})$ the left-hand side of $(\ref{3-3})$ goes to zero when $|E|\to 0$, uniformly with respect to $\ge_n$.
This implies that $(\ref{ineq3})$ is still valid and the conclusion of the proof is as in \rth{regN=2}.
\qeda
\medskip

\noindent\Remark A simpler statement which implies $(\ref{3-1})$ is the following.
\begin{equation}\label{3-4}
\lim_{\gd\to 0}\myint{0}{\gd}\left(\myint{B_r(z)}{}V(y)(\gf(y)-\ge)_+^2dy\right)\myfrac{dr}{r^{N+1}}=0,
\end{equation}
uniformly with respect to $0<\ge\leq \ge_0$ and to $z\in\Gs_\ge$. The proof is similar to the one of \rprop{CS}.
\medskip

\noindent\Remark When the function $V$ depends essentially of the distance to $\prt\Gw$ in the sense that
\begin{equation}\label{ineq3'}\BA {l}
\abs {V(x)}\leq v(\gf(x))\qquad\forall x\in\Gw,
\EA\end{equation} 
and $v$ satisfies
\begin{equation}\label{ineq4}\BA {l}
\myint{0}{a}tv(t)dt<\infty,
\EA\end{equation} 
Marcus and V\'eron proved \cite[Lemma 7.4]{MV3} that $\CR(u)=\prt\Gw$, for any positive solution  $u$ of $(\ref{F1})$.  This assumption implies also $(\ref{3-1})$.  The proof is similar to the one of \rprop{Tr}.
%%%%%%%%%%%%%%%%%%%%%%%%%%%%%%%%%%%%%%%%%%%%%%%%%%%%%%%%%%%%%%%%%%%%%%%%%%%%%%%%%%%%%%%%%%%%%%%%%%%%%%%%%%%%%%%%%%%%%%%%%%%%SUB-SECTION%%%%%%%%%%%%%%%%%%%%%%%%%%%%%%%%%%%%%%%%%%%%%%%%%%%%%%%%%%%%%%%%%%%%%%%%%%%%%%%%%%%%%%

\subsection{The sweeping method}
 
This method introduced in \cite{RV} for analyzing isolated singularities of solutions of semilinear equations has been adapted in  \cite{MV0} and \cite{MV4} for defining an extended trace of positive solutions of differential inequalities in particular in the super-critical case. Since the boundary  trace of a positive solutions of $(\ref {F1})$ is known on $\CR(u)$ we shall study the sweeping with measure concentrated on the singular set $\CS(u)$

\bprop{inf} Let $u\in C(\Gw)$ be a positive solution of $(\ref {F1})$ with singular boundary set $\CS(u)$. If  $\gm\in\frak M_+(\CS(u))$ we denote $v_\gm=\inf\{u,u_\gm\}$. Then
 \begin{equation}\label{F2'}
-\Gd v_\gm+V(x)v_\gm\geq 0\qquad\text{in }\Gw,
\end{equation}
and $v_\gm$ admits a boundary trace $\gg_u(\gm)\in \frak M_+(\CS(u))$. The mapping $\gm\mapsto \gg_u(\gm)$ is nondecreasing and $\gg_u(\gm)\leq\gm$.
\es
\Proof By \cite{sta}, $(\ref{F2'})$ holds.
But $Vu_\gm\in L^1_\gf(\Gw)\Longrightarrow Vv_\gm\in L^1_\gf(\Gw)$, if we set $w:=\BBG[Vv_\gm]$, then
$v_\gm+w$ is nonnegative and super-harmonic, thus it admits a boundary trace in $\frak M_+(\prt\Gw) $ that we denote by $\gg_u(\gm)$. Clearly $\gg_u(\gm)\leq \gm$ since $v_\gm\leq u_\gm$  and $\gg_u(\gm)$ is nondeacreasing with $\gm$ as $\mu\mapsto u_\gm$ is. Finally, since $v_\gm$ is a supersolution, it is larger that the solution of $(\ref{F1})$ with the same boundary trace $\gg_u(\gm)$, and there holds
 \begin{equation}\label{F3}
u_{\gg_u(\gm)}\leq v_\gm.
\end{equation}

\bprop{Borel} Let
 \begin{equation}\label{F4}
\gn_{_S}(u):=\sup\{\gg_u(\gm):\gm\in   \frak M_+(\CS(u))\}.
\end{equation}
Then $\gn_{_S}(u)$ is a Borel measure on $\CS(u)$. 
\es
\Proof We borrow the proof to Marcus-V\'eron \cite{MV4}, and we naturally extend any positive Radon measure to a positive bounded and regular Borel measure by using the same notation. It is clear that $\gn_{_S}(u):=\gn_{_S}$ is an outer measure in the sense 
 that 
 \begin {eqnarray}\label {borel subadd}\gn_{_S}(\emptyset)=0,\;\mbox { and }\gn_{_S} (A)\leq\sum_{k=1}^\infty\nu (A_{k}),\;\mbox { whenever 
 }A\subset\bigcup_{k=1}^\infty A_{k}.\end {eqnarray}
 Let $A$ and $B\subset \CS(u)$ be disjoint Borel subsets. In 
order to prove that
\begin {eqnarray}\label {borel add}
\gn_{_S} (A\cup B)=\gn_{_S} (A)+\gn_{_S} (B),
\end {eqnarray}
we first notice that the relation holds if $\max \{\gn_{_S} (A),\gn_{_S} 
(B)\}=\infty$. Therefore we assume that $\gn_{_S} (A)$ and $\gn_{_S} (B)$ are 
finite. For $\varepsilon>0$ there exist two bounded positive
measures $\mu_{1}$ and $\mu_{2}$ such that 
$$\gamma_{u}(\mu_{1})(A)\leq \nu (A)\leq \gamma_{u}(\mu_{1})(A)+\varepsilon/2
$$
and
$$\gamma_{u}(\mu_{2})(B)\leq \nu (B)\leq \gamma_{u}(\mu_{2})(B)+\varepsilon/2
$$
Hence
$$\begin {array}{l}
\gn_{_S}(A)+\gn_{_S} (B)\leq \gamma_{u}(\mu_{1})(A)+
\gamma_{u}(\mu_{2})(B)+\varepsilon\\
\phantom {\gn_{_S}(A)+\gn_{_S} (B)}
\leq\gamma_{u}(\mu_{1}+\mu_{2})(A)+\gamma_{u}(\mu_{1}+\mu_{2})(B)
+\varepsilon\\
\phantom {\gn_{_S}(A)+\gn_{_S} (B)}
= \gamma_{u}(\mu_{1}+\mu_{2})(A\cup B)+\varepsilon\\
\phantom {\gn_{_S}(A)+\gn_{_S} (B)}
\leq \gn_{_S} (A\cup B)+\varepsilon.
\end {array}
$$
Therefore $\gn_{_S}$ is a finitely additive measure. If $\{A_{k}\}$ ($k\in 
\BBN$) is a sequence of of disjoint Borel sets and $A=\cup A_{k}$, then
$$\gn_{_S} (A)\geq \gn_{_S} \left(\bigcup_{1\leq k\leq n}A_{k}\right)=\sum_{k=1}^n\gn_{_S} 
(A_{k})\Longrightarrow \gn_{_S} (A)\geq \sum_{k=1}^\infty\gn_{_S} (A_{k}).
$$
By $(\ref {borel subadd})$, it implies that $\gn_{_S}$ is a countably additive measure. \qeda

\bdef{extend} The Borel measure $\gn(u)$ defined by 
 \begin{equation}\label{ext-tr}
\gn(u)(A):=\gn_{_S}(A\cap\CS(u))+\gm_u(A\cap\CR(u)),\qq\forall A\subset\prt\Gw,\, A\text{ Borel},
\end{equation}
is called the extended boundary trace of $u$, denoted by $Tr^{e}(u)$.
\es

\bprop {loc}If $A\subset \CS(u)$ is a Borel set, then
 \begin{equation}\label{rest}
\gn_{_S}(A):=\sup\{\gg_u(\gm)(A):\gm\in  \frak M_+(A)\}.
\end{equation}
\es
\Proof If $\gl,\gl'\in\mathfrak M_+(\CS(u))$  
$$\inf\{u,u_{\gl+\gl'}\}=\inf\{u,u_{\gl}+u_{\gl'}\}\leq \inf\{u,u_{\gl}\}+\inf\{u,u_{\gl'}\}.$$
Since the three above functions admit a boundary trace, it follows that
$$\gg_u(\gl+\gl')\leq \gg_u(\gl)+\gg_u(\gl').
$$
If $A$ is a Borel subset of $\CS(u)$, then $\gm=\gm_{A}+\gm_{A^c}$  where $\gm_{A}=\chi_{_E}\gm$. Thus
$$\gg_u(\gm)\leq \gg_u(\gm_{A})+\gg_u(\gm_{A^c}),
$$
and
$$\gg_u(\gm)(A)\leq \gg_u(\gm_{A})(A)+\gg_u(\gm_{A^c})(A).
$$
Since $\gg_u(\gm_{A^c})\leq \gm_{A^c}$ and $\gm_{A^c}(A)=0$, it follows
$$\gg_u(\gm)(A)\leq \gg_u(\gm_{A})(A).
$$
But $\gm_{A}\leq \gm$, thus $\gg_u(\gm_{A})\leq \gg_u(\gm)$ and finally
 \begin{equation}\label{rest1}
\gg_u(\gm)(A)= \gg_u(\gm_{A})(A).
\end{equation}
If $\gm\in  \frak M_+(A)$, $\gm=\gm_{A}$, thus $(\ref{rest})$ follows.\qeda

\bprop{vanish} There always holds 
 \begin{equation}\label{F5}
\gn(u)({\mathcal S}{\scriptstyle ing}_{_ V}(\Omega ))=0,
\end{equation}
where ${\mathcal S}{\scriptstyle ing}_{_ V}(\Omega )$ is the vanishing set of $K_V^\Gw(x,.)$ defined by $(\ref{van1})$.
\es
\Proof This follows from the fact that for any $\gm\in\frak M_+(\prt\Gw)$ concentrated on 
${ \mathcal S}{\scriptstyle ing}_V(\Omega )$, $u_\gm=0$. Thus $\gg_u(\gm)=0$. If $\gm$ is a general measure, we can write $\gm=\chi_{_{{\mathcal S}{\scriptstyle ing}_{_ V}(\Omega )}}\gm+
\chi_{_{({\mathcal S}{\scriptstyle ing}_{_ V}(\Omega ))^c}}\gm$, thus $u_\gm=u_{\chi_{_{({\mathcal S}{\scriptstyle ing}_{_ V}(\Omega ))^c}}\gm}$. Because of $(\ref{F3})$ 
$$\gamma_u(\gm)({\mathcal S}{\scriptstyle ing}_{_ V}(\Omega ))=\gamma_u(\chi_{_{({\mathcal S}{\scriptstyle ing}_{_ V}(\Omega ))^c}}\gm)({\mathcal S}{\scriptstyle ing}_{_ V}(\Omega ))\leq (\chi_{_{({\mathcal S}{\scriptstyle ing}_{_ V}(\Omega ))^c}}\gm)({\mathcal S}{\scriptstyle ing}_{_ V}(\Omega ))=0,
$$
thus $(\ref{F5})$ holds.\qeda\medskip

\noindent {\Remark} This process for determining the boundary trace is ineffective if there exist positive solutions $u$ in $\Gw$ such that
$$\lim_{\gd_{\Gw}(x)\to 0}u(x)=\infty.
$$
This is the case if $\Gw=B_R$ and $V(x)=c(R-\abs x)^{-2}$ ($c>0$). In this case $K^\Gw_V(x,.)\equiv 0$. For any $a>0$, there exists a radial solution of
\begin{equation}\label{F6'}
-\Gd u+\myfrac{cu}{(R-|x|)^2}=0\qquad \text {in } B_R
\end{equation}
under the form 
 \begin{equation}\label{F6}
u(r)=u_a(r)=a+c\myint{0}{r}s^{1-N}\myint{0}{s}u(t)\myfrac {t^{N-1}dt}{(R-t)^2}.
\end{equation}
Such a solution is easily obtained by  fixed point, $u(0)=a$ and the above formula shows that $u_a$ blows up when $r\uparrow R$. We do not know if there a exist non-radial positive solutions of $(\ref{F6'})$.
More generaly, if $\Gw$ is a smooth bounded domain, we do not know if there exists a non trivial positive solution of
\begin{equation}\label{F6''}
-\Gd u+\myfrac{c}{d^{2}(x)}u=0\qquad \text {in } \Gw.
\end{equation}

\bth{regN<2} Assume $V\geq 0$ and satisfies $(\ref{cond'})$. If $u$ is a positive solution of $(\ref{F1})$, then  $Tr^{e}(u)=\gn(u)$ is a bounded measure.\es
%%%%%%%
\Proof Set $\gn=\gn(u)$ and assume $\gn(\prt\Gw)=\infty$. By dichotomy there exists a decreasing sequence of relatively open domains $D_n\subset\prt\Gw$ such that $\overline D_n\subset D_{n-1}$, diam$\,D_n=r_n\to 0$ as $n\to\infty$, and $\gn(D_n)=\infty$. For each $n$, there exists a Radon measure $\gm_n\in \frak M_+(D_n)$ such that  $\gg_u(\gm_n)(D_n)=n$, 
and 
$$u\geq v_{\gm_n}=\inf\{u,u_{\gm_n}\}\geq u_{\gg_u(\gm_n)}.
$$
Set $m_n=n^{-1}\gg_u(\gm_n)$, then $m_n\in \frak M_+(D_n)$ has total mass $1$ and it converges in the weak*-topology to $\gd_a$, where $\{a\}=\cap_nD_n$. By \rth{stab}, $u_{m_n}$ converges to 
$u_{\gd_a}$. Since $u\geq nu_{m_n}$, it follows that 
$$u\geq \lim_{n\to\infty}nu_{m_n}=\infty,$$
a contradiction. Thus $\gn$ is a bounded Borel measure (and thus outer regular) and it corresponds to a unique Radon measure.\qeda
\medskip

\noindent\Remark If $N=2$, it follows from \rth{regN=2} that $u=u_\gn$ and thus the extended boundary trace coincides with the usual boundary trace. The same property holds if $N\geq 3$, if $(\ref{3-1})$ holds.
%%%%%%%%%%%%%%%%%%%%%%%%%%%%%%%%%%%%%%%%%%%%%%%%%%%%%%%%%%%%%%%%%%%%%%%%%%%%%%%%%%%%%%%%%%%%%%%%%%%%%%%%%%%%%%%%%%%%%%%%%%%%%%%%%%%%%%%%%%%%%%%%%%%%%%%%%%%%%%%%%APPENDIX%%%%%%%%%%%%%%%%%%%%%%%%%%%%%%%%%%%%%%%%%%%%%%%%%%%%%%%%%%%%%%%%%%%%%%%%%%%%%%%%%%%%%%%%%%%%%%%%%%%%%%%%%%%%%%%%%%%%%%%%%%%%%%%%%%%%%%
\appendix\section{Appendix: A necessary condition for the fine regularity of a boundary point with respect to a Schr\"odinger equation  } 
\setcounter{equation}{0}
{\large by Alano Ancona}\footnote{D\'{e}partement de Math\'{e}matiques, B\^atiment 425, Universit\'{e} Paris-Sud 11, Orsay 91\hspace{0.5pt}405 France \\
 Email address: alano.ancona@math.u-psud.fr}

{\bf \large }

\noindent \hspace{6mm}
\noindent \hspace{6mm}

This appendix is devoted to the derivation of   a sufficient condition  --stated in  \rth{thmain} below (section A1)-- for the {\em fine singularity} of a boundary point of a Lipschitz domain with respect to a potential $V$. This theorem answers a question  communicated  by Moshe Marcus and Laurent V\'{e}ron to the author --and related to the work \cite{MV1} by Marcus and V\'{e}ron-. The expounded proof goes back to the unpublished manuscript \cite{anc4}. In a forthcoming paper other criterions for fine regularity will be given -- in particular  a simple  explicit  necessary and sufficient condition for the fine regularity of a boundary point and a criteria for having almost everywhere regularity in a subset of the boundary.

The exposition  can be read independently  of the above paper of L.\ V\'{e}ron and C.\ Yarur. The few  notions necessary to the statement of \rth{thmain} are recalled in section A1. Section A2 is devoted to some known basic  preliminary results and the proof of  \rth{thmain} is given in section A3.

{\sl Acknowledgment.} The author is grateful to Moshe Marcus and Laurent V\'{e}ron for  bringing to his attention  their motivating question.

\subsection{ Framework,  notations and  main result} 

Let $\Omega $ be a bounded Lipschitz domain in ${ \mathbb R} ^N$. Denote  $ \delta_{\Omega } (x):=d(x;{ \mathbb R} ^N\setminus \Omega )$ 
the distance from $x$ to the complement of $\Omega$ in ${ \mathbb R} ^N$ and for $a>0$, let ${ \mathcal V}(\Omega ,a)$ denote the set of all nonnegative measurable function $V:\Omega \to  { \mathbb R} $ such that $V(x) \leq a/(\delta _{\Omega }(x))^2$ in $\Omega $. We also let $x_0$ to denote a fixed reference point in $\Omega $. 

For  $V \in   { \mathcal V}(a,\Omega )$, we will consider the Schr\" odinger  operator $L_V:= \Delta  -V$ associated with the potential $V$. Here $ \Delta  $ is the classical Laplacian in ${ \mathbb  R}^N$.

\noindent {\sl The kernels $K_y$, $ \tilde K_y^V$ and $K_y^V$.}
It is  well known (\cite{HW1}, \cite{HW2}) that  to each point $y \in   \partial \Omega $ corresponds a unique positive harmonic function $K_y$ in $\Omega $ that vanishes on $\partial \Omega $ and satisfies the normalization condition $K_y(x_0)=1$. This function is the Martin kernel w.r.\ to the Laplacian  in $\Omega $  with pole at $y$ and normalized at $x_0$. It may  also be seen as a Poisson kernel with respect to $ \Delta  $ in $\Omega $.

\noindent 
The function $K_y$ is obviously superharmonic in $\Omega $ with respect to $L_V$ and we may hence consider its  greatest $L_V$-harmonic minorant $ \tilde K_y^V$ in $\Omega $ defining hence another kernel function  at $y$.

By the results in \cite{anc3} (see paragraph {\bf A2}  below)  it is also known that for each $y \in   \partial \Omega $ there exists a unique positive  $L_V$-harmonic function $K_y^V$ in $\Omega $ that vanishes on $\partial \Omega \setminus  \{ y \}$ and satisfies $K_y^V(x_0)=1$. Thus $ \tilde K_y^V=c_y\, K_y^V$ with $c_y= \tilde K_y^V(x_0)$. Here a function $u: \Omega \to  { \mathbb  R}$ is $L_V$-harmonic if $u$ is the continuous representative of a weak solution $u$ of $L_V(u)=0$ (so $u \in   H^1_{loc}(\Omega )$ by assumption and necessarily $u \in   W^{2,p}_{loc}(\Omega )$ for all $p< \infty $).

The set of ``{\em finely}" regular boundary points   with respect to $L_V $ in $\Omega $ is 
  \begin{align}{ \mathcal R}{\scriptstyle eg}_{_ V}(\Omega )&:= \{ y \in   \partial \Omega\,;\,  \tilde K_y^V>0 \}= \{ y \in   \partial \Omega \,;\, c_y>0\,  \}   \end{align}
\noindent -since $c$ is u.s.c.\ this is a $K_ \sigma  $ subset of $\partial \Omega $- and the set of ``{\em finely}"  irregular boundary points is ${\mathcal S}{\scriptstyle ing}_{_ V}(\Omega ):= \partial \Omega\setminus { \mathcal R}{\scriptstyle eg}_{_ V}(\Omega )$. These notions were introduced by E.\ B.\ Dynkin in his study  of  positive solutions in $\Omega $ of a non linear equation such as $ \Delta  u= u^q $, $q >1$ -in which case, given $u$,  we recover Dynkin's definition on taking $V= \vert  u \vert  ^{ q -1}$. See the books \cite{Dbook1}, \cite{dyn2} of E.\,B.\ Dynkin and the references there. From the probabilistic point of view, a boundary point $y \in   \partial \Omega$ is $L_V $ finely regular iff for the Brownian motion  $ \{  \xi  _s \}_{0 \leq s< \tau }$ starting say at $x_0$ and conditioned to exit from $\Omega$ at  $y$, it holds that $\int _0^ \tau \, V( \xi  _s)\, ds <+ \infty $ a.s., or in other words, iff the probability for this  process to reach $y$ when killed at the rate $e^{-V( \xi  _s)\, ds}$   is strictly positive.

Let us now state  \rth{thmain}. It  answers  the question (2005) of Marcus-V\'{e}ron  alluded to above: 
suppose that for sufficiently many Lipschitz path (resp.\ every linear path)  $ \gamma :[0, \eta  ] \to  \overline  \Omega $ such that $ \gamma (0)=y$ and $d( \gamma (t),\partial \Omega ) \geq c\,   \vert  \gamma (t)-y \vert $ for $0 \leq t \leq  \eta  $ and some $c>0$, it holds that 
$$\int _0^ \eta  \, t\, V( \gamma (t))\, dt =+ \infty; $$
 does it follow that $y$ is finely singular w.r.\ to $V$ and $\Omega $ ?

\bth{thmain} Let $y \in   \partial \Omega$ and let $C_{\ge,y}:= \{ x \in   \Omega\,;\,  \delta_\Omega  (x) \geq  \varepsilon \, d(x,y) \}$ for $\; 0<\varepsilon <1$.
If  \begin{equation} \int _{ C_{\ge,y} }\, V(x)  \,{ \frac {  dx} { \vert  x-y \vert  ^{N-2} }}=+ \infty \end{equation}
for some $ \varepsilon >0 $,
then $y \in   { \mathcal S}{\scriptstyle ing}_{_V}(\Omega )$.\es

\subsection{Boundary Harnack principle for $L_V$}

To prove  \rth{thmain} we will rely on the main result of \cite{anc3} (see also \cite{ancbook}) in well-known forms more or less explicit    in \cite{anc3} (see e.g.\ Theorem 5$'$ and Corollary 27  there) or \cite{ancbook}. In this section we state these needed  ancillary results and fix some notations to be used in  what follows.

  Fix positive reals $r,\,  \rho >0$ such that $0<10\,r< \rho $ and let $f$ be a ${ \frac {   \rho } {10r }} $ lipschitz function in the ball $B_{N-1}(0,r)$ of ${ \mathbb  R}^{N-1}$ -- we let $B_{N-1}(m,s)$ to denote the ball in ${ \mathbb  R}^{N-1}$ of center $m$ and radius $r$--. Define then the region $U_f(r, \rho )$ in ${ \mathbb  R}^N$ as follows 
\begin{equation} U_f(r, \rho ):= \{ (x',x_N) \in   { \mathbb  R}^{N-1}\times { \mathbb  R}\simeq { \mathbb  R}^N\,;\,  \vert  x' \vert  <r,\, f(x')<x_N< \rho \,  \}    \end{equation}
\noindent We will also denote it $U$ (leaving $f$, $r$ and  $ \rho $ implicit) when convenient. Set   $\partial _{\#}U:=\partial U \cap  \{ x=(x',x_N) \in   { \mathbb  R}^N\,;\,  \vert  x' \vert   \leq r,\, x_N=f(x')\,  \}$ and define $T(t):=B_{N-1}(0;tr)\times (-t \rho , +t \rho )$ .

Recall ${ \mathcal V}_a(U)$ is the set of all Borel nonnegative functions $V$ in $U$
such that $V(x) \leq  { \frac {  a} { \delta (x)^2 }}$ for  $x \in   U$. For $V$ H\"older continuous (in fact for a natural class of second order elliptic operators) the following statement goes back to \cite{anc1}. See also \cite {Dah} for $V=0$.

\blemma{lemmaBHP}  Let $V \in   { \mathcal V}_a(U)$ and set $L _V:=  \Delta  -V$.  There is a constant $C$ depending only on $N$, $a$ and ${ \frac {   \rho } {r }}$ such that for  any  two positive $L_V $-harmonic functions  $u$ and $v$ in $U $ that vanish on $\partial_\# U $, 
 \begin{align} { \hspace {10truemm}  \frac {u(x)} {u(A) }} &\leq  C\; { \frac {  v(x)} {v(A) }} \hspace {5truemm} \ \ {for\  all\ } x \in   U \cap T({\frac  { 1} {2}}) 
\end{align}
where $A=A_U=(0,\dots, 0,{ \frac {   \rho } {2 }})$.
\es

\noindent {   \em Proof.} Let us briefly recall -for readers convenience- how this lemma follows from Theorem 1 in \cite{anc3}.
By homogeneity we may assume that $r=1$ and that $ \rho $ is fixed.  Let $A'=(0,\dots, 0,{ \frac {  2 \rho } {3 }})$ and let $B_N$ denote the  open ball $B_N(0,1)$ in ${ \mathbb  R}^N$. 
It is easy to construct  a bi-Lipschitz map  $F: U \to  B_N(0,1)$ with a bi-lipschitz constant depending only on $ \rho $ and $N$ and which   maps $A'$ onto $0$,  $U \cap T(1/2)$ onto  $B_N^-:= \{ x \in   B_N\,;\, x_N<-{ \frac { 1  } {2}} \}$ and $U\setminus T({\frac  { 3} {4}} )$ onto $B_N^+:= \{ x \in   B_N\,;\, x_N \geq { \frac { 1  } {2}}\,  \}$. 

Standard calculations show  that if $u$ is $ \Delta  -V$ harmonic in $U$ then the function $u_1:=u \circ F^{-1}$ is ${L}_1-V \circ F^ {-1}$ harmonic in $B_N$  for some (symmetric) divergence form elliptic operator $L_1= \sum _{i,j} \partial _i(a_{ij} \partial _j)$  in $B_N$ satisfying $ C_1^{-1}\; I_N \leq   \{ a_{ij} \}  \leq C_1\, I_N$ with $C_1=C_1(N,{\frac  { r} { \rho }} ) \geq 1$. Let $V_1=V \circ F_1^{-1}$. Clearly $V_1 \in   { \mathcal V}(B_N,a')$ for  $a'=C^2\, a$.

Other simple calculations show that the operator ${ \mathcal L}=(1- \vert  x \vert ) ^2(L_1-V_1)$ seen as a map $H^1_{loc}(B_N) \to  H^{-1}_{loc}(B_N)$ is an adapted elliptic operator in divergence form over the hyperbolic ball $B_N$ (i.e.\ w.r.\ to the hyperbolic metric $ds^2= {\frac  {  \vert  dx \vert  ^2} {(1- \vert  x \vert^2  )^2}} $) in the sense of \cite{anc3}. Moreover since the form $ \varphi   \mapsto  \int_{B_N}  a_{ij}\partial _i \varphi \, \partial _j \varphi  \, dx - \varepsilon _0\int_{B_N}  {\frac  {   \varphi  ^2} {(1- \vert  x \vert  )^2}} \, dx$ is coercive for $ \varepsilon _0= \varepsilon _0(C_1,N)>0$ chosen sufficiently  small, the differential  operator ${ \mathcal L}$ is weakly coercive which means that there exists $ \varepsilon _0= \varepsilon _0(N,{\frac  { r} { \rho }})>0$ such that ${ \mathcal L}+ \varepsilon _0$ admits a Green's function in $B_N$.

This shows that  Theorem 1 in \cite{anc3} applies to ${ \mathcal L}$. Thus there is a constant $c=c( \varepsilon_0 ,C_1,N)$, $c \geq 1$, such that for $z=(z',z_N)  \in    B_N ^+$  and $y \in   B_N^-$ one has  \begin{align}\label{eqappendixendix1} c^{-1}\, G_{ \mathcal L}(y,z) &\leq G_{ \mathcal L}(y,0)\,G_{ \mathcal L}(0,z) \leq c\, G_{ \mathcal L}(y,z)  
\end{align} Here we have also used the standard  Harnack inequalities for ${ \mathcal L}$ and have denoted $G_{ \mathcal L}$ the ${ \mathcal L}$ Green's function in $B_N$ w.r. to the hyperbolic metric (we adopt the notational convention that $u(x):=G_{ \mathcal L }(x,y)$ satisfies ${ \mathcal L}u=- \delta _x$ in the weak sense \cite{sta} w.r.\ to the hyperbolic volume). Notice that $G_{ \mathcal L}(x,y)=   \delta (y)^{N-2} g(x,y)$ if $g$ is Green's function of $ L_1-V_1$ in $B_N$ (w.r.\ to the usual metric).

Supppose that  $u_1$ is positive ${ \mathcal L}$ harmonic (i.e.\  $L_1-V_1$ harmonic) in $B_N$ and that $u_1$ vanishes on $\partial B_N  \cap     \{ x \in  \partial B_N\,;\, x_N \leq {\frac  { 1} {2}}  \}$. Then  $u_1$ can be represented as a Green potential in $B_N \cap     \{x\,;\, x_N<{\frac  { 1} {2}}   \}$ : $ u_1(y)= \int  G_{ \mathcal L} (y,z)\, d\nu (z)$ where $\nu $ is a nonnegative Borel measure on $ \{ z \in   B_N\,;\, z_N={\frac  { 1} {2}}  \}$ and $y_N \leq {\frac  { 1} {2}} $. So upon integrating (\ref{eqappendixendix1}) we get (with another constant $c$)   \begin{align}   
 \label{eqappendix2} c^{-1}\, u_1(y) &\leq u_1(0)\,g(y,0) \leq c\, u_1(y)  
\end{align} for $y \in  B_N^-$.  Thus if $u$ is a positive $L_V$ solution in $U$ that vanishes in $\partial _\#U$ it follows --on using  the change of variable $y=F(x)$-- that 
\begin{align}   
 \label{eqappendix3} c^{-1}\, u(x) &\leq u(A')\,G(x,A') \leq c\, u(x)  
\end{align} for $x  \in   U({\frac  { 1} {2}} )$, 
 where $G$ is Green's function w.r.\ to $L_V$ in $U$. Using Harnack inequalities for $L_V$, the lemma easily follows.$\square$\medskip

\noindent\Remark {\rm Using \rlemma{lemmaBHP}, well known  arguments (see \cite{anc1})  show that for every bounded Lipschitz  domain  $\Omega $ in ${ \mathbb R} ^N$ and every $V \in   { \mathcal V}(\Omega ,a)$, $a>0$, the following potential theoretic properties hold in  $\Omega $ equipped with  $L_V:= \Delta  -V$ (we let $G_y^V$ to denote the $L_V$ Green's function in $\Omega $ with pole at $y$) 
: \hspace{5truemm} (a) For each $P  \in   \partial \Omega $, the limit $K_P^V(x)=\lim_{y  \to    P} G^V_y(x)/G^V_y(x_0)$, $x \in   \Omega $, exists and $K_P^L$ is a positive $L_V$-harmonic function $K_P^L$  in $\Omega $ which depends continuously on $P$ and vanishes continuously in $\partial \Omega \setminus  \{P \}$, 
\hspace{5truemm} (b) For each $P  \in   \partial \Omega $, every positive $L_V$-solution in $\Omega $ that vanishes on $\partial \Omega \setminus  \{ P \}$ is  proportional to $K_P^V$,  \hspace{5truemm} (c) Every positive $L_V$-solution $u$  in $\Omega $ can be written in a unique way as $u(x)=\int _{\partial \Omega }\, K_P^V(x)\, d\mu  (P)$, $x \in   \Omega $, for some positive (finite) measure $\mu  $ in $\partial \Omega  $. See \cite{anc3}.}

\subsection{Proof of \rth{thmain}}

Again $\Omega $ is a bounded Lipschitz domain in ${ \mathbb  R}^N$  and $V \in   { \mathcal V}(\Omega ,a)$, $a \geq 0$.

 For the proof we use a simple variant of the comparison principle given in  \rlemma{lemmaBHP}. Notations are as before, in particular   $U=U_f(r, \rho )$ is the domain considered in {\bf A2} and $A=A_U=(0,\dots, 0,{ \frac {   \rho } { 2}})$. Let $A'= (0,\dots, 0,{ \frac {   2\rho } { 3}})$.

\blemma {vBHI} Let $u$ be  positive harmonic (w.r.\ to $ \Delta  $) in $U$, let $v$ be positive $ \Delta  -V$-harmonic in $U$ and assume  that $u=v=0$ in $\partial_\# U $. Then 
\begin{align} { \frac {v(x)} {v(A) }} &\leq c\; { \frac {  u(x)} {u(A) }} \hspace {5truemm} \ \ {for\  \ } x \in   U \cap T({\frac  { 1} {2}}) \end{align}
\flushleft for some positive constant $c$ depending only on $ \rho /r$, the constant $a$ and $N$.
\es

\noindent {\sl Proof.} We have seen that $v(x) \leq c\,v(A')\, G_{A'}^V(x)$ in $U \cap T({ \frac { 1  } {2}})$ and we know that $G_A^V \leq G_A^0$ in $U$ if $G_{A'}^V$ is $( \Delta  -V)$-Green's function in $U$ with pole at $A'$. By maximum principle, Harnack inequalities and the known behavior of $G_{A'}^0$ in $B(A', { \frac { r  } {4}})$ (more precisely $G_A^0(x) \leq c_1:=c_1(r, N)$ in $\partial B(A',{ \frac { r  } {4}})$) we have that $u(x) \geq c_1\,v(A)\, G_{A'}^0(x)$ in $U\setminus B(A',{ \frac { r  } {4}})$. So that --using  Harnack inequalities in $B(A', { \frac { r  } {2}})$ for $u$ and $v$-- the lemma follows. $\square$\medskip

\noindent \Remark {\rm The opposite estimate, i.e.  $ { \frac {u(x)} {u(A) }} \leq  C\; { \frac {  v(x)} {v(A) }}$ (with another constant $C>0$), cannot be expected to hold in general as shown by simple (and obvious) examples}.

\noindent Denote $g_{x_0}^V $ the Green's function with respect to $ \Delta  -V$ in $\Omega $ and with  pole at $x_0$. For $y \in   \partial \Omega $, a pseudo-normal for $\Omega $ at $y$  is a unit vector $\nu  \in { \mathbb  R}^N$ such that that for some small $  \eta >0$, the set $C(y,\nu _y,  \eta ):= \{ y+t(\nu _y+v)\, ;\, 0<t<  \eta ,\,  \Vert v \Vert  \leq     \eta  \; \} $ is contained in $   \Omega $.  

\bprop{appenequivdef} Given $y \in   \partial \Omega$ and a pseudo-normal $\nu _y$ at $y$ for $U$, the following assertions are equivalent: \smallskip

\noindent (i) $\tilde K^V_y=0$ (i.e.\, $y \in   { \mathcal S}{\scriptstyle ing}_V(\Omega )$) \smallskip

\noindent (ii) $\limsup_{t\downarrow 0} K_y^V(y+t\nu _y)/K_y(y+t\nu _y)=+ \infty $ \smallskip

\noindent (iii) $\lim_{t\downarrow 0} K_y^V(y+t\nu _y)/K_y(y+t\nu _y)=+ \infty $ \smallskip

\noindent (iv) $\lim_{x \to  y} g_{x_0}^V(x)/g^0_{x_0}(x)=0$. 

\es
\Proof 
(a) We  first recall a standard consequence of \rlemma{lemmaBHP} that relates $g^V_{x_0}$ and $K_y^V$ near $y$ (for any $y \in   \partial \Omega $). 

Consider $u=K^V_y $ and $v:=g^V_{y+t\nu_y }$. Using \rlemma{lemmaBHP} and the fact that $v\sim t^{2-N}$ in $\partial B(y+t\nu_y ,  {\frac  {  \eta  } {2}}  t)$, $0<t< \eta  $,    we see that $u(x) \sim u(y+t\nu_y )\, t^{N-2}\, g_{y+t\nu_y }^V(x)$ for $x \in   \Omega \setminus B(y+t\nu _y, t  \eta /2   )$ (here $\sim$ means ``is in between two constant times " with constants depending only on $y$, $\Omega $, $\nu_y $ and $a$). 

Taking in particular $x=x_0$ we obtain that $K^V_y(y+t\nu_y )\sim  1  /(t^{N-2}g^V(y+t\nu_y ;x_0)) $. In particular considering the special case $V=0$, we get also that $K_y(y+t\nu_y )\sim  1  /(t^{N-2}g(y+t\nu_y ;x_0)) $.

(b) Using the above we see that (ii) is equivalent to (iv)$'$: $\liminf_{t\downarrow 0} g_{x_0}^V(y+t\nu _y)/g_{x_0}^0(y+t\nu _y)=0 $. 

(c) Now to show that $(iv)$ and $(iv)'$ are equivalent  we may assume that $y=0$, $\nu _y=(0,\dots, 0,1)$ and  (with the notations above in {\bf A.2}) that $T(1) \cap \Omega =U$, $U=U_f(r, \rho )$ and $x_0 \in   \Omega \setminus \overline  U$.

 Applying \rlemma{vBHI} to $U$,  $u=g^V_{x_0}$,$v=g_{x_0}$, and $U_t=U_{t_j}$ for a sequence $t_j$, $t_j\downarrow 0$ such that  $u(A_{t_j})=o(v(A_{t_j}))$, $A_{t_j}=(0,\dots, 0,t_j)$,  we get that $u(x) \leq c{\frac  { u(A_{t_j})} {v(A_{t_j})}}\, v(x)$ in $\Omega  \cap    T(t_j{\frac  {  \rho } {2}} )$. Hence (iv)$'$ imply (iv).  And  --using (a) again-- conditions (ii), (iii) and (iv) are equivalent.

(d) Similarly if on the contrary $g^V(A_j, x_0) \geq \, c\, g(A_j, x_0)$, for some sequence $A_j=t_j\nu $, $t_j\downarrow 0$ and a positive real $c$, we have (since a priori $g^V \leq g$) that  : 
\begin{equation} K_{A_j}^V(x):=g^V(A_j,x)/g^V(A_j,x_0)  \leq c^{-1}\, K_{A_j} (x)=c^{-1}\,g(A_j,x)/g(A_j,x_0)  
\end{equation}
and letting $j \to   \infty $ we get $K_y^V \leq c^{-1} K_y$. Thus, (i)$\Rightarrow$(iv).

Since obviously (ii) $\Rightarrow$ (i), \rprop{appenequivdef} is proved.   $\square$

The  next lemma is the key for the proof of  \rth{thmain}. Returning again to the canonical Lipschitz domain $U=U_f(r, \rho )$, let $V \in   { \mathcal V}_a(U)$ and for $ \theta   \in   (0,{ \frac {  1} {10 }})$, let $U ^\theta  := \{ x \in   U\,;\, d(x,\partial U) \geq  \theta  r\, \}$, $I_U^ \theta  :=\int _{U^ \theta  } \, V(x)\, { \frac {  dx\ \ } {  \vert  \delta _U(x) \vert  ^{N-2} }}$. 

Obviously ${ \frac {  1} {r^{N-2}  }} \int _{U^ \theta  } V(x)\, dx \leq I_U^ \theta    \leq \, { \frac {  1} { (\theta  r)^{N-2}  }} \int _{U^ \theta  } V(x)\, dx$.

\blemma{keyappenlemma} Let $u $, $ \tilde u $ be two nonnegative  continuous functions in $\overline U$ that are respectively $ \Delta  $-harmonic and $L_V$-harmonic in $U$. Assume  that $ \tilde u \leq u$ in $\partial U$ and $ \tilde u=u=0$ in $\partial_\# U$. Then for some constant $c=c( { \frac { r  } { \rho }}, a, \theta  , N)>0$,
\begin{equation} (1+cI_ \theta  )\,  \tilde u(x)\,  \leq u (x) \hspace{3truemm} {\text for \  } x \in   U \cap T({ \frac {  1} {2 }}) \end{equation} 
\es
\Proof Since the assumptions and the conclusion are  invariant under dilations we may assume that $r$ is fixed as well as $ \rho $. Replacing $u$ by the harmonic function in $U$ with same boundary values as $ \tilde u$ we may also assume that $u= \tilde u$ in $\partial U$. Since $ \Delta  (u- \tilde u)=-V\,  \tilde u$ and $ u-\tilde u$ vanishes on $\partial U$, we see that $u- \tilde u=G_U(V  \tilde u)$ where $G_U$ is the usual Green's function in $U$.

By Harnack property and since $G_U(x,y) \geq c=c( \theta  ,a,N)>0$ for $x \in   B_1=B(A_1,{ \frac { r  } {100}})$, $A_1=(0,\dots, 0,{\frac  { 3r} {4}}  )$, and $y \in   U^\theta  $, we have
\begin{align} u(x)- \tilde u(x) \geq c\,I_ \theta  \,    \tilde u(A_1), \ x  \in   B_1.  \nonumber
\end{align}
Thus in $U$, $w(x):=u(x)- \tilde u(x) \geq c\,\, I_ \theta  \, \,  \tilde u(A_1)\, R_1^{B_1}(x)$ where $R_1^{B_1}$ is the (classical) capacitary potential (\cite{doob}) of $B_1$ in $U$ and  using the comparison principle Lemma 1 for $V=0$ we have $w \geq c\, I_ \theta  \,  \tilde u(A_1)\,\, {\frac  {u } {u(A_1)}} $ in $U({ \frac {  1} { 2}}):=T({ \frac {  1} {2 }}) \cap U$. 

Using then \rlemma{lemmaBHP} (and Harnack inequalities)
\begin{align} w(x) \geq c''\,  I_ \theta  \;  \tilde u(A_1)\,\, {\frac  {  \tilde u(x)} { \tilde u(A_1)}} \;= c'''\, I_ \theta  \,  \tilde u(x),\,\; x \in U({\frac  { 1} {2}}) \nonumber   
\end{align}
Thus, $u(x)  \geq (1+c'''\, I_ \theta  )\, \tilde u(x) $ in $U({\frac  { 1} {2}} )$. $\square$

\vspace{1truecm}
\noindent {\bf Proof of \rth{thmain}.} We may assume that $y=0$, that for some $r$, $ \rho $, $f$, $\Omega  \cap    T(1)=U:=U_f(r, \rho )$ (with the notation fixed above in section {\bf A2}) and that $x_0\notin \overline U$. 

Set $T_n=T(2^{-n})$, $C_y^n:=  C_{\ge,y} \cap    (T_n\setminus T_{n+1})$ for $n \geq 1$, $u=G_{x_0}^0$, $ \tilde u=G_{x_0}^V$ (where $G_{x_0}^V$ is Green's function with pole at $x_0$ with respect to $ \Delta  -V $ in $\Omega $).
One may also observe  that  $ \varepsilon $ may be assumed so small that $ \Sigma _ 0^\varepsilon $ contains the truncated cone $C:= \{ (x',x_N)\,;\, x_N<{ \frac {   \rho } {2 }},\,  \vert  x' \vert  <{ \frac {   r} { \rho  }}\,x_N\,  \}$. 

For each $n \geq 0$ there is a greatest $ \alpha _n>0$ such that $u \geq  \alpha _n\,  \tilde u$ in $U_n$ (we know that $ \alpha _n \leq 1$). By the key \rlemma{keyappenlemma} (and elementary geometric considerations)
\begin{align}  \alpha _{n+1}  \geq  \alpha _n\, (1+cI_{n+1})   \;\; {\rm if \ } I_m:=\int _{C_m} {\frac  { V(x)} {\, \delta _\Omega (x)^{N-2}}}\; dx  
\end{align}
for some constant $c=c(  \varepsilon   , {\frac  { r} { \rho }},a,N) $ independent of $n$. Thus 
\begin{align}  \alpha _n \geq  \alpha _0\,  \prod_{k=1}^{n-1} (1+c\, I_k)  \geq  \alpha _0\,(1+c\,  \sum _{k=1} ^{n-1} I_k) \geq c\, \alpha _0\, \int _{C_1\setminus C_{n+1}}{\frac  { V(x)} { \delta _\Omega (x) ^{N-2}}} \, dx \nonumber
\end{align}
which shows that $\lim  \alpha _n=+ \infty $. Thus $G_{x_0} ^V=o(G_{x_0}^0)$ at $y$ and by \rprop{appenequivdef} the point $y$ belongs to ${ \mathcal S}ing_V(\Omega )$. $\square$\medskip

%%%%%%%%%%%%%%%%%%%%%%%%%%%%%%%%%%%%%%%%%%%%%%%%%%%%%
%%%%%%%%%%%%%%%%%%%%%%%%%%%%%%%%%%%%%%%%%%%%%%%%%%%%%
%%%%%%%%%%%%%%%%%%%%%%%%%%%%%%%%%%%%%%%%%%%%%%%%%%%%%
%%%%%%%%%%%%%%%%REFERENCES%%%%%%%%%%%%%%%%%%%%%%
%%%%%%%%%%%%%%%%%%%%%%%%%%%%%%%%%%%%%%%%%%%%%%%%%%%%%
%%%%%%%%%%%%%%%%%%%%%%%%%%%%%%%%%%%%%%%%%%%%%%%%%%%%%
%%%%%%%%%%%%%%%%%%%%%%%%%%%%%%%%%%%%%%%%%%%%%%%%%%%%%
%%%%%%%%%%%%%%%%%%%%%%%%%%%%%%%%%%%%%%%%%%%%%%%%%%%%%

%%%%%%%%%%%%%%%%%%%%%%%%%%%%%%%%%%%
%%%%%%%%%%%%%%%%%%%%%%%%%%%%%%%%%%%
%%%%%%BEGINING OF THE ARTICLE%%%%%%%%%%%%%%
%%%%%%%%INTRODUCTION%%%%%%%%%%%%%%%%%%%%%%%%%%%%%%%%%%%%%%%%%%%%%%%%%%%%%
%%%%%%%%%%%%%%%%%%%%%%%%%%%%%%%%%%%

%Furthermore%Furthermore%Furthermore
%Furthermore%Furthermore%Furthermore
%Furthermore%Furthermore%Furthermore
 %%END DOCUMENT%%%%%%%%%%%%%%%%%%%%%%%%
\end {document}